\documentclass[journal]{IEEEtran}
\usepackage{url}
\usepackage{cite}
\usepackage{latexsym}
\usepackage{amsmath}
\usepackage{amstext,amsfonts,epsf,epsfig,pifont}
\usepackage{graphics,bm,amsmath}
\newlength{\intwidth}

\def\XXint#1#2#3{{\setbox0=\hbox{$#1{#2#3}{\int}$}
\vcenter{\hbox{$#2#3$}}\kern-.5\wd0}}

\newtheorem{theorem}{Theorem}

\newtheorem{definition}{Definition}[section]

\newtheorem{assumption}{Assumption}[section]

\newcommand{\ftp}{\Longleftrightarrow}

\newcommand{\seq}{\!\!\!=\!\!\!}

\newcommand{\calE}{\ensuremath{\mathcal{E}}}

\newcommand{\calR}{\ensuremath{\mathcal{R}}}

\begin{document}

\noindent \emph{The following statements are placed here in accordance with the copyright policy of the Institute of Electrical and Electronics Engineers, Inc., available online at}
\url{http://www.ieee.org/web/publications/rights/policies.html}.\\

\noindent
Lilly, J. M., \&  Olhede, S. C. (2010). On the analytic wavelet \indent transform. \emph{IEEE Transactions on Information Theory,} \textbf{56} \\ \indent (8), 4135--4156.\\

\noindent This is a preprint version.  The definitive version is available from the IEEE or from the first author's web site, \url{http://www.jmlilly.net}.\\

\noindent \copyright 2010 IEEE. Personal use of this material is permitted. However, permission to reprint/republish this material for advertising or promotional purposes or for creating new collective works for resale or redistribution to servers or lists, or to reuse any copyrighted component of this work in other works must be obtained from the IEEE.\\

\newpage

\title{On the Analytic Wavelet Transform}
\author{Jonathan~M.~Lilly,~\IEEEmembership{Member,~IEEE,}
and Sofia~C.~Olhede,~\IEEEmembership{Member,~IEEE}
\thanks{The work of J. M. Lilly was supported  by award \#0526297 from the Physical Oceanography program of the United States National Science Foundation. A collaboration visit by S.~C.~Olhede to Earth and Space Research in the summer of 2006 was funded by the Imperial College Trust.}
\thanks{J.~M.~Lilly is with Earth and Space Research, 2101 Fourth Ave., Suite 1310, Seattle, WA 98121, USA (e-mail: lilly@esr.org).}
\thanks{S.~C.~Olhede is with the Department of Statistical Science, University College London, Gower Street,
London WC1 E6BT, UK (e-mail: s.olhede@ucl.ac.uk).}}

\markboth{IEEE Transactions on Information Theory}{Lilly \& Olhede}

\maketitle
\begin{abstract}
An exact and general expression for the analytic wavelet transform of a real-valued signal is constructed, resolving the time-dependent effects of non-negligible amplitude and frequency modulation.  The analytic signal is first  locally represented as a modulated oscillation, demodulated by its own instantaneous frequency, and then Taylor-expanded at each point in time.  The terms in this expansion, called the instantaneous modulation functions, are time-varying functions which quantify, at increasingly higher orders, the local departures of the signal from a uniform sinusoidal oscillation.  Closed-form expressions for these functions are found in terms of Bell polynomials and derivatives of the signal's instantaneous frequency and bandwidth. The analytic wavelet transform is shown to depend upon the interaction between the signal's instantaneous modulation functions and frequency-domain derivatives of the wavelet, inducing a hierarchy of departures of the transform away from a perfect representation of the signal.  The form of these deviation terms suggests a set of conditions for matching the wavelet properties to suit the variability of the signal, in which case our expressions simplify considerably.  One may then quantify the time-varying bias associated with signal estimation via wavelet ridge analysis, and choose wavelets to minimize this bias.
\end{abstract}
\IEEEpeerreviewmaketitle

\section{Introduction}
\PARstart{T}{his} paper derives properties of the Analytic Wavelet Transform (AWT), a special family of complex-valued wavelet transforms, for the analysis of modulated oscillations. The complex-valued wavelet transform has emerged as an important non-stationary signal processing tool; a discussion of its properties together with references may be found in \cite{selesnick05-ispm}.  Continuous complex wavelets have been used for the characterization of modulated oscillatory signals \cite{delprat92-itit,mallat,carmona97-itsp,carmona99-itsp,scheper03-itsp} and discontinuities \cite{tu05-itit}. The applications of complex wavelets to analysis of real signals includes mechanical vibratory signals \cite{kim05-jsv,zhang03-itec}, seismic signals \cite{olhede05-icae}, position time series from drifting oceanic buoys \cite{lilly06-npg}, and quadrature Doppler signals in blood flow \cite{aydin00-ispl}. Much attention has also focused on the design of discrete wavelet filters that approximate the effect of an analytic continuous wavelet transform, with  important contributions by Kingsbury \cite{kingsbury01-acha}, Selesnick \cite{selesnick02-itsp}, and others \cite{gopinath03-itsp,fernandes05-itip}. In addition, the signal form we discuss in this article has been used to model speech \cite{mcaulay86-itassp} and echolocation of bats \cite{simmons71-science} as well as  gravitational waves \cite{anderson99-prd}. The results regarding the properties of the AWT derived in this article will thus be relevant to the analysis of signals from a number of fields.

A broad class of interesting signals may be modeled as modulated oscillations, with the analytic signal as the foundation \cite{gabor46-piee,vakman77-spu}.  One wishes to recover the properties of the signal, without specifying a parametric model for its structure, based on a generally noisy or contaminated observation.  Of particular interest are the time-varying amplitude and phase of the analytic signal as well as their derivatives.  For contaminated signals the direct construction of the analytic signal via the Hilbert transform can lead to disastrous results as the amplitude and phase will then reflect the aggregate properties of the multi-component signal \cite{loughlin01-itsp}.  It is necessary to isolate the signal of interest while simultaneously rendering it analytic.   The AWT is a recipe for constructing a family of versions of a time series which are both localized and analytic.

A second analysis step, termed wavelet ridge analysis \cite{delprat92-itit,mallat}, then identifies a special set of points from which the properties of an underlying analytic signal can be accurately estimated. This method can exhibit excellent performance, principally on account of its insensitivity to signal contamination due to the time/frequency localization of the wavelets. But the price of the localization is that the analytic signal is no longer precisely recovered in the absence of contamination, unlike direct construction of the analytic signal---i.e., bias is introduced in the estimation procedure. The departure of the estimated analytic signal from the true analytic signal is negligible if the signal modulation over the time support of the wavelet is also negligible \cite{mallat}, a strong constraint since many real-world signals exhibit substantial modulation.

The purpose of this paper is to determine the exact properties of the AWT, and the resulting ridge-based signal estimation, for local analysis of oscillations with non-negligible modulation. Although Mallat~\cite{mallat} derives error bounds for analytic signal estimation in the weakly modulated case, time-dependent errors due to moderate or strong modulation have not yet been considered.  Understanding this modulation-induced bias is important in order to correctly interpret the amplitude and frequency estimates provided by the wavelet ridge analysis.  These results may also be used as a guide in choosing wavelets which explicitly minimize bias effects.

The structure of the paper is as follows. Necessary background is presented in Section~\ref{section:background}. Section~\ref{localsection} then introduces a novel representation of a modulated oscillatory signal as a series of departures from a pure sinusoidal oscillation.  In Section~\ref{transformsection}, a general expression for the AWT of a modulated oscillation is then derived.  This result is used in Section~\ref{ridgesection} to examine the deterministic bias properties of ridge-based signal estimation.  The paper concludes with a discussion.

All numerical code associated with this paper is made freely available for use by others, as noted in Appendix~\ref{section:software}.

\section{Background}\label{section:background}
This section reviews the specification of the amplitude and phase of a modulated oscillatory signal via the analytic signal, and their estimation by the wavelet ridge method \cite{delprat92-itit,mallat} using a general family of analytic wavelets  \cite{olhede02-itsp,lilly09-itsp}.

\subsection{Modulated Oscillations}

A real-valued amplitude- and frequency-modulated signal may be usefully represented as \cite{voelcker66a-ieee}
\begin{eqnarray}
x(t)&=& a_+(t) \cos \phi_+(t)
\label {modulatedmodel}
\end{eqnarray}
with the amplitude $a_+(t)$ and phase $\phi_+(t)$ defined in terms of the \emph{analytic signal} \cite{gabor46-piee}.  The analytic signal is specified in the frequency domain by
\begin{equation}
x_+(t)\equiv 
\label{analyticfilter}
\frac{1}{2 \pi}\int_{-\infty}^\infty X_+(\omega)e^{i\omega t}\,d\omega
\end{equation}
where we have introduced
\begin{equation}
X_+(\omega)\equiv 2U(\omega) X(\omega)
\end{equation}
with $U(\omega)$ being the Heaviside unit step function.

The construction of the analytic signal $x_+(t)$ permits the amplitude $a_+(t)$  and  phase $\phi_+(t)$ to be uniquely defined\footnote{Note that at isolated points where $a_x(t)=0$, the value of the phase is typically defined by continuity \cite{picinbono97-itsp}.} via
\begin{equation}
a_+(t) e^ {i\phi_+(t)}\equiv x_+(t)
\end{equation}
and the original signal is recovered by $x(t)= \Re\left\{x_+(t)\right\}$.  While more than one amplitude and phase pair may yield the same real-valued signal in (\ref{modulatedmodel}), the symbols $a_+(t)$ and $\phi_+(t)$ denote the so-called \emph{canonical amplitude} and \emph{canonical phase}\cite{picinbono97-itsp} associated with the analytic signal.  The conditions under which a given amplitude and phase can be recovered in this manner have been examined by \cite{bedrosian63-ire}.  The rates of change of the amplitude and phase are quantified by
\begin{eqnarray}
\label{instfreq}
\upsilon(t)&\equiv &\Re\left\{ \frac{d}{dt}\,\ln x_+(t)\right\}=\frac{a_+'(t)}{a_+(t)}\\
\omega(t)&\equiv &\Im\left\{ \frac{d}{dt}\,\ln x_+(t)\right\}=\phi_+'(t)
\end{eqnarray}
which are referred to as the \emph{instantaneous bandwidth} \cite{cohen} and \emph{instantaneous frequency} \cite{voelcker66a-ieee,voelcker66b-ieee,boashash92a-ieee}, respectively.  These two fundamental instantaneous quantities have an intimate connection to the first two frequency-domain moments of the signal's spectrum; see e.g. \cite{lilly10a-itsp} and references therein.

It frequently arises that one wishes to estimate the instantaneous properties---$a_+(t)$, $\phi_+(t)$, $\upsilon(t)$ and $\omega(t)$---of a modulated oscillatory signal $x(t)$ believed to be present in a noisy observation.  Typically one is presented with an observed time series $x^{(o)}[t_n]$ at discrete times $t_n\in t_1$, $t_2$, $\ldots, t_N$
\begin{equation}
x^{(o)}[t_n]=x(t_n)+x^{(\epsilon)}[t_n]
\end{equation}
where $x(t_n)$ is a discretely sampled modulated oscillation and  $x^{(\epsilon)}[t_n]$ is a discrete noise process.  Given the noisy observed signal $x^{(o)}[t_n]$, one wishes to estimate the properties of the modulated oscillation  $x(t)$.  A powerful method for accomplishing this task is \emph{wavelet ridge analysis}, described subsequently.  This method extracts an estimate of a modulated oscillation from the analytic wavelet transform.  In the estimation procedure, there are three sources of error:
\begin{enumerate}
\item[(i)] Errors associated with the discrete sampling;
\item[(ii)] Random errors due to the noise $x^{(\epsilon)}[t_n]$; and
\item[(iii)] Errors dependent upon the oscillation  $x(t)$ itself.
\end{enumerate}
The purpose of this paper is to examine this third type of error, which may be called ``bias''. Henceforth we assume that the noise process $x^{(\epsilon)}[t_n]$ vanishes and that the sampling is perfect, and consequently we work in continuous time.

\subsection{The Analytic Wavelet Transform}\label{definitionssection}
A wavelet  $\psi(t)\in L^2({\mathbb{R}})$ is an analyzing function used to localize a signal simultaneously in time and frequency.  By definition, a wavelet has zero mean and finite energy, and additionally satisfies the ``admissibility condition'' \cite{holschneider}
\begin{eqnarray}
\int_{-\infty}^\infty \frac{\left|\Psi(\omega)\right|^ 2}{\left|\omega\right|}\, d\omega&< & \infty
\end{eqnarray}
where $\Psi(\omega)$ is the Fourier transform of the wavelet.   The wavelet transform of a signal $x(t)\in L^2({\mathbb{R}})$ is a series of projections onto rescaled and translated versions of  $\psi(t)$
\begin{equation}
W_\psi(t,s)\equiv \int_{-\infty}^{\infty} \frac{1}{s}\, \psi^*\left( \frac{\tau -t}{s}\right) x(\tau)\,d\tau \label{wavetrans}
\end{equation}
which are indexed by both the time parameter $t$ and a scale parameter $s$;  the asterisk denotes the complex conjugate. Note the choice of a $1/s$ normalization rather than the more common $1/\sqrt{s}$, as we find the former to be more convenient for oscillatory signals.

Here we will consider only wavelets that vanish for negative frequencies, i.e. that have $\Psi(\omega)=0$ for $\omega<0$. Such wavelets are called \emph{analytic}\footnote{This terminology reflects the fact that, if $\Psi(\omega)$ has no support on negative frequencies, $\psi(z)$ will be an analytic function of a complex argument~$z$; see the discussion in Appendix~1 of \cite{olhede04-biometrika}.} and  (\ref{wavetrans}) then defines the \emph{analytic wavelet transform} (AWT).  Alternatively the AWT may be represented in the frequency domain as
\begin{equation}
W_\psi(t,s)
=\frac{1}{2\pi}\int_{0}^{\infty} \Psi^*(s\omega) X(\omega)\,e^{i\omega t}\,d\omega \label{wavetransfourier}
\end{equation}
with the integration requiring only positive frequencies on account of the exact analyticity of the wavelet.  The analytic wavelet $\Psi(\omega)$ has a maximum amplitude in the frequency domain at $\omega=\omega_\psi$, which is called the ``peak frequency''.  We choose to set the value of the wavelet at the peak frequency to $\Psi(\omega_\psi)= 2$, since then with $x(t)=a_o \cos\left(\omega_o t\right)$ we have the convenient result $\left|W_\psi\left(t,\omega_\psi/\omega_o\right)\right| =\left|a_o\right|$.

A useful way to categorize wavelet behavior is through normalized versions of the derivatives of the frequency-domain wavelet.  We define the wavelet's \emph{dimensionless derivatives} as
\begin{equation}
\widetilde\Psi_{n}(\omega) \equiv  \omega^n\frac{\Psi^{(n)}(\omega)}{\Psi(\omega)}\label{normalizedwavelet}
\end{equation}
where the superscript ``$(n)$'' denotes the $n$th derivative.  We will consider only wavelets for which the second derivative at the peak frequency  $\Psi''(\omega_\psi)$ is real-valued, in which case it is also negative since $\Psi(\omega_\psi)$ is a maximum by definition.  Then
\begin{equation}
P_\psi\equiv\sqrt{-\widetilde\Psi_{2}(\omega_\psi)} =\sqrt{ -\omega_\psi^2\frac{\Psi''(\omega_\psi)}{\Psi(\omega_\psi)}}\label{pdef}
\end{equation}
defines a real-valued quantity $P_\psi$ which is a nondimensional measure of the wavelet \emph{duration}.

$P_\psi/\pi$ may be shown to be the number of oscillations at the peak frequency $\omega_\psi$ which fit within the central time window of the wavelet \cite{lilly09-itsp}, as measured by the standard deviation of the demodulated wavelet $\psi(t)e^{-i\omega_\psi t }$.  Approximating the wavelet by its second-order Taylor expansion about $\omega_\psi$ gives
\begin{multline}
\Psi(\omega)\approx\Psi(\omega_\psi)+\frac{1}{2}\left(\omega-\omega_\psi\right)^2\Psi''(\omega_\psi)
\\=\Psi(\omega_\psi)\left[1-\frac{1}{2}\left(\frac{\omega}{\omega_\psi}-1\right)^2 P_\psi^2\right],\quad\omega\approx\omega_\psi\label{pbandwidth}
\end{multline}
and so we have $\Psi\left(\omega_\psi(1\pm 1/P_\psi)\right)/\Psi(\omega_\psi)\approx 1/2$.  Thus the inverse duration $1/P_\psi$ can also be seen as a nondimensional measure of the wavelet bandwidth.  After this second-order description of the wavelet,  $\widetilde\Psi_{3}(\omega_\psi)$ offers the next-higher-order description, and can be interpreted as quantifying the asymmetry of the wavelet about its peak frequency \cite{lilly09-itsp}.

\subsection{A General Family of Analytic Wavelets}\label{section:gmw}
To examine the role of the analyzing wavelet in shaping the performance of wavelet ridge analysis, we will need a general family of analytic wavelets whose properties may readily calculated.  The generalized Morse wavelets  \cite{daubechies88-ip,olhede02-itsp,lilly09-itsp} are such a family.  These wavelets are defined in the frequency domain by
\begin{eqnarray}
\psi_{\beta,\gamma}(t)\ftp\Psi_{\beta,\gamma}(\omega)&=& U(\omega) \,a_{\beta,\gamma}\, \omega^\beta e^{-\omega^\gamma}
\label{morse}
\end{eqnarray}
where $a_{\beta,\gamma}$ is a normalizing constant and $U(\omega)$ is again the unit step function.    The generalized Morse wavelets are controlled by two parameters, $\beta$ and $\gamma$, the roles of which in shaping wavelet properties were examined in detail by~\cite{lilly09-itsp}.  To be a valid wavelet one must have $\beta>0$ and $\gamma>0$.  By varying these two parameters, the generalized Morse wavelets can be given a broad range of characteristics while remaining exactly analytic. Note that we will replace the subscript ``$\psi$'' with ``$\beta,\gamma$'' to denote quantities pertaining to these wavelets.

Simple expressions for important properties of the generalized Morse wavelets are given in \cite{lilly09-itsp}. The peak frequency is $\omega_{\beta,\gamma} \equiv\left(  \beta/\gamma \right) ^{1 /\gamma}$, and choosing $a_{\beta,\gamma}\equiv  2 (e\gamma/\beta)^{\beta/\gamma}$, with $e$ being Euler's number $2.7182\ldots$, then gives $\Psi_{\beta,\gamma}(\omega_{\beta,\gamma})=2$.  The dimensionless duration $P_{\beta,\gamma}$ of the generalized Morse wavelets is
\begin{equation}
P_{\beta,\gamma}\equiv\sqrt{-\widetilde\Psi_{2;\beta,\gamma}(\omega_{\beta,\gamma})} =\sqrt{\beta\gamma}
\end{equation}
while the third-order dimensionless derivative at the peak frequency is
\begin{equation}
\widetilde\Psi_{3;\beta,\gamma}(\omega_{\beta,\gamma})=
-(\gamma-3)P_{\beta,\gamma}^2 \label{morsepsi3}.
\end{equation}
A general expression for $\widetilde\Psi_{n;\beta,\gamma}(\omega_{\beta,\gamma})$ of any order may be found in \cite{lilly09-itsp}.

In fact, these wavelets form a very broad family that subsumes many other types of wavelets.  It was shown by \cite{lilly09-itsp} that the generalized Morse wavelets encompass two other popular families of analytic wavelets: the Cauchy or Klauder wavelet family ($\gamma=1$) and the analytic ``Derivative of Gaussian'' wavelets ($\gamma=2$).   The diversity of the generalized Morse wavelets is due to the fact that they are a function of two parameters, $\beta$ and $\gamma$, hence their second-order and third-order properties $P_{\beta\gamma}$ and $\widetilde\Psi_{3;\beta,\gamma}(\omega_{\beta,\gamma})$ may be independently varied.

\begin{figure*}[t!]
        \noindent\begin{center}\includegraphics[height=7in,angle=-90]{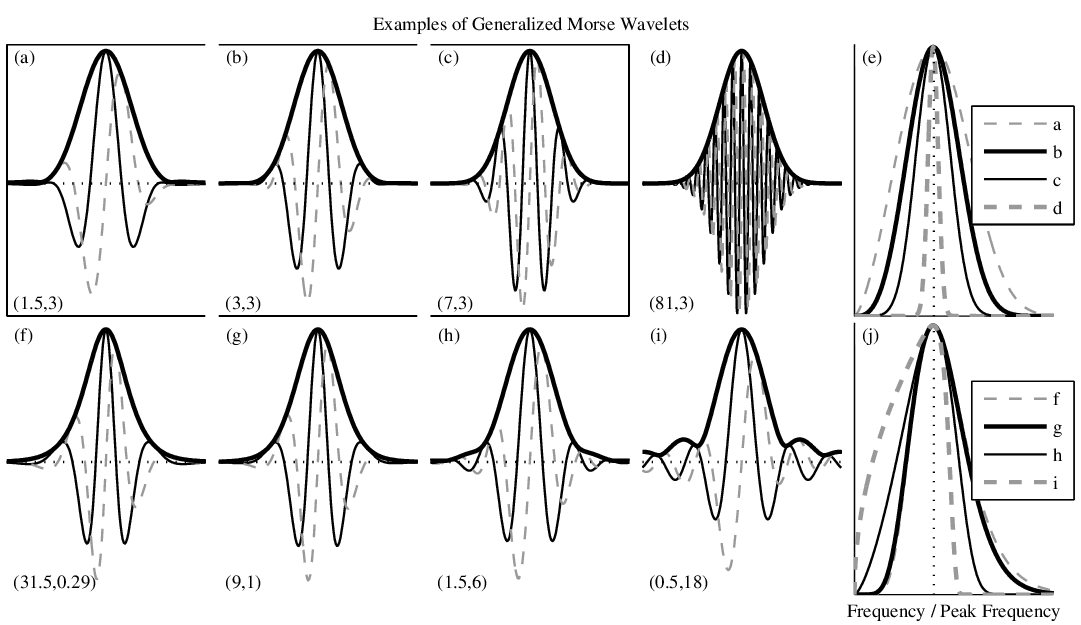}\end{center}
        \caption{\footnotesize
Examples of the generalized Morse wavelets.  Panels (a--d) and (f--i) show wavelets in the time domain for different $(\beta,\gamma)$ values, which are indicated in the lower left corner of each panel; for presentation, the wavelets are rescaled by their maximum amplitude.  The real part is a solid line, the imaginary part is dashed, and the modulus is a thick solid line.  The frequency-domain versions of the wavelets in the top and bottom rows are then given in panels (e) and (j) respectively.  Three wavelet which will be used in the future, (a--c), are set off with a box.  The upper row of wavelets shows the effect of increasing $\beta$ with $\gamma$ fixed at $\gamma=3$, hence $\widetilde\Psi_{3;\beta,\gamma}(\omega_{\beta,\gamma})$ vanishing; here $P_{\beta,\gamma}=\sqrt{\beta\gamma}$ takes on values of $\sqrt{4.5}$, $\sqrt{9}=3$, $\sqrt{21}$, and $\sqrt{243}$. The lower row shows the effect of increasing $\gamma$ and decreasing $\beta$ such that $P_{\beta,\gamma}$ remains fixed at $P_{\beta,\gamma}=3$;
here $\widetilde\Psi_{3;\beta,\gamma}(\omega_{\beta,\gamma})/\widetilde\Psi_{2;\beta,\gamma}(\omega_{\beta,\gamma})$ is equal to -2.7, -2.0, 3.0, and  14.9, respectively.}\label{analytic-morsies}
\end{figure*}

Examples of the generalized Morse wavelets are shown in Fig.~\ref{analytic-morsies}.  The upper row shows the effect of increasing $\beta$ with $\gamma$ fixed at $\gamma=3$, hence vanishing $\widetilde\Psi_{3;\beta,\gamma}(\omega_{\beta,\gamma})$, with $P_{\beta,\gamma}=\sqrt{\beta\gamma}$ increases from left to right.  The wavelet becomes more oscillatory in the time domain, or more tightly peaked in the frequency domain.  The lower row shows the effect of increasing $\gamma$ and decreasing $\beta$  with $\widetilde\Psi_{3;\beta,\gamma}(\omega_{\beta,\gamma})$  increasing from negative values on the left to positive values on the right.  With $P_{\beta,\gamma}$ fixed, the number of oscillations within the central window does not change, but the long-time behavior of the wavelet changes considerably; in the frequency domain, an enhancement to the right of the peak shifts to the left of the peak as $\gamma$ increases.  This illustrates the effect of separately varying second-order and third-order wavelet properties.

\normalsize

\subsection{Wavelet Ridge Analysis}\label{section:ridgeintro}

The idea of wavelet ridge analysis \cite{delprat92-itit,mallat} is that there exist special time/scale curves, called \emph{wavelet ridge curves} or simply \emph{ridges}, along which properties of a modulated oscillatory signal are accurately represented. Unlike the instantaneous frequency curve, which is not known, the ridge curves are based on properties of the transform itself and can be readily located. The AWT evaluated along the ridge constitutes an estimator for a presumed modulated oscillatory signal.

A ridge curve is based on an aggregation of points known as \emph{ridge points}.  Two separate definitions of ridge points are in common use.

\begin{definition}{Ridge Points}\\
An {\em amplitude ridge point} of $W_\psi\left(t,s\right)$ is a time/scale pair $\left(t,s\right)$ satisfying the two conditions
\begin{eqnarray}
\frac{\partial}{\partial s}\, \Re\left\{\ln W_\psi\left( t , s\right)\right\}& = &  0\label{ampridge}\\
\frac{\partial^2}{\partial s^ 2}\, \Re\left\{\ln W_\psi\left( t , s\right)\right\}& < &  0.
\label{ampridge2}
\end{eqnarray}
Since $\Re\left\{\ln W_\psi\left( t , s\right)\right\} = \ln \left|W_\psi\left( t , s\right)\right|$, these conditions state that for each fixed time point $t$, an amplitude ridge point corresponds to the scale at which a local maximum in the transform magnitude occurs.

Similarly a {\em phase ridge point} of $W_\psi\left(t , s\right)$ is a time/scale pair $\left(t,s\right)$ satisfying the two conditions
\begin{eqnarray}
\frac{\partial}{\partial t}\, \Im\left\{\ln \label{phaseridge} W_\psi\left(t,s\right)\right\}-\frac{\omega_\psi}{s}& = &  0\\
\frac{\partial}{\partial s}\left[\frac{\partial}{\partial t}\, \Im\left\{\ln W_\psi\left(t,s\right)\right\}- \frac{\omega_\psi}{s}\right]& < & 0.\label{phaseridge2}
\end{eqnarray}
Condition (\ref{phaseridge}) states that the rate of change of transform phase matches $\omega_\psi/s$, which is interpretable as a frequency associated with scale $s$; see  \cite{lilly09-itsp} for a discussion of this point.  Note that condition (\ref{phaseridge2}), like (\ref{ampridge2}), has the effect of identifying amplitude maximum rather than minima.  While (\ref{phaseridge2}) is not standard, we introduce it here in order that the phase ridges may be defined without reference to the transform amplitude.
\end{definition}

Ridge points are then grouped into sets called {\em ridge curves}.

\begin{definition}{Ridge Curves}\\
Let the set of all amplitude ridge points of some real-valued signal $x(t)$ with respect to a wavelet $\psi(t)$ be denoted $S^{\{a\}}$, while $S^{\{p\}}$ denotes the set of all phase ridge points.  Henceforth we will use the notation such as $S^{\{\cdot\}}$, with a superscript ``$\cdot$'' referring to either ``$a$'' or ``$p$''.  Then a {\em ridge curve} $s^{\{\cdot\}}(t)$ is a scale curve as a function of time which maps out a contiguous collection of individual ridge points.  The ridge curve is defined over some time interval $T^{\{\cdot\}}$  and is constrained to additionally satisfy the continuity condition
\begin{equation}
\left|\frac{d}{dt}\,s^{\{\cdot\}}(t)\right|<\infty.
\end{equation}
This latter condition excludes discontinuities in the scale curve $s^{\{\cdot\}}(t)$ as a function of time, as well as multiple values of scale at a particular time.
\end{definition}

The union of all ridge points is also known as the wavelet \emph{skeleton} of the signal [\citen{antoine}, p14--18].  An estimate of a modulated oscillatory signal may then be constructed by evaluating the wavelet transform along the ridge curve.  Here we have assumed the presence of a single modulated oscillation; signals consisting of a superposition of such oscillations may be treated in a similar fashion provided the instantaneous frequency curves are sufficiently separated in time and/or frequency \cite{carmona99-itsp}.

\begin{definition}{The Ridge-Based Signal Estimate}\\
The amplitude or phase {\em ridge-based signal estimate} is
\begin{equation}
\widehat x^{\{\cdot\}}_{+,\psi}(t)\equiv W_\psi\left( t, s^{\{\cdot\}}(t)\right)\quad t\in  T^{\,\{\cdot\}}\label{estimator}
\end{equation}
which is the set of values the wavelet transform takes along the ridge curve.  This simple form is due to the $1/s$  normalization and the choice $\Psi(\omega_\psi) = 2$ introduced in Section~\ref{definitionssection}.
\end{definition}

It was shown by \cite{mallat} that the error in $\widehat x^{\{\cdot\}}_{+,\psi}(t)$ becomes negligible when $\upsilon(t)$, $\upsilon'(t)$, and $\omega'(t)$, together with a fourth term involving broadband bias, all tend to zero.  However, the time-varying form of the error for non-vanishing modulation strength, and the conditions governing an appropriate choice of wavelet for a given signal, have not yet been examined.

\subsection{Application to Oceanographic Data}

\begin{figure*}[t!]
        \noindent\begin{center}\includegraphics[width=7in,angle=0]{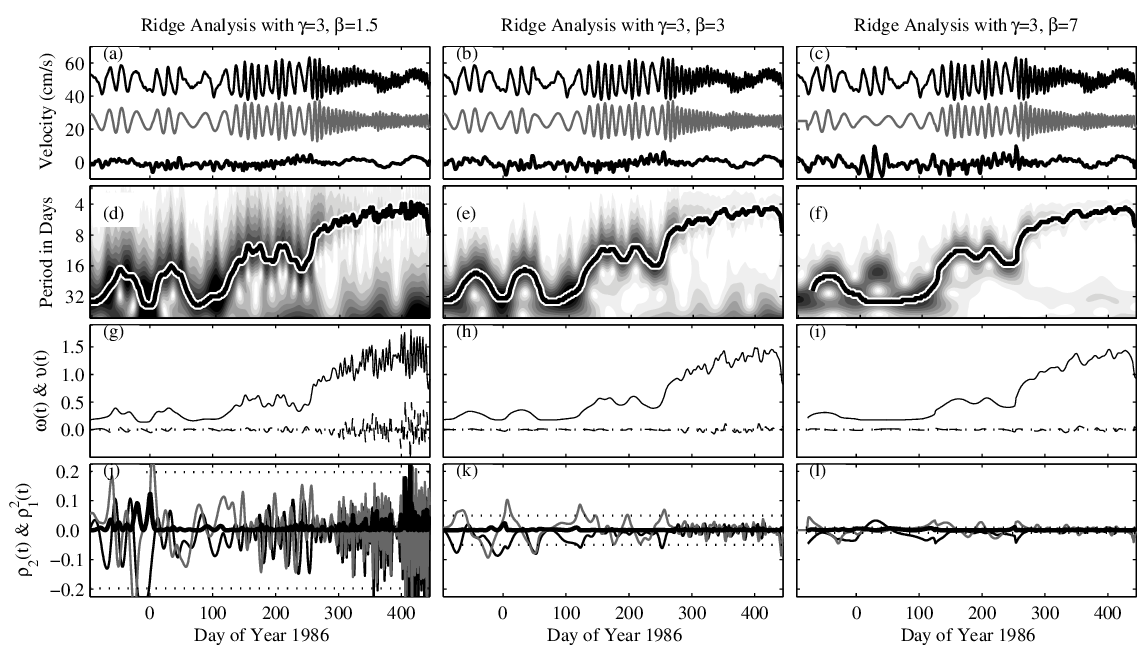}\end{center}
        \caption{\footnotesize
        Example of wavelet ridge analysis.  The first row (a--c) shows the data together with the resulting ridge-based signal estimate and the residual; these are offset with the data (which is the same for all panels) at the top, the estimated signal in gray in the middle, and the residual at the bottom.  The three columns correspond to the use of three different wavelets, specifically those shown in Fig.~\ref{analytic-morsies}a--c.  The three wavelet transforms are shown in (d--f) with the ridge curves marked as heavy lines.   The final two rows are used in a subsequent section.  The third row (g--i) shows the instantaneous frequency (solid) and bandwidth (dashed) of the estimated signals; these quantities are estimated as described later in Section~\ref{ridgesection}.  The fourth row (j--l) shows the estimated real (solid) and imaginary (dashed) parts of the second instantaneous modulation function $\widetilde\rho_2(t)$, along with the contribution $\upsilon^2(t)/\omega^2(t)$ due to the squared bandwidth (heavy curve).  Dotted lines are plotted at $\pm4/P_{\beta,\gamma}^4$, which as discussed later is a measure of the degree of variability appropriate for each wavelet.
}\label{analytic-transforms}
\end{figure*}

An example of wavelet ridge analysis is shown in Fig.~\ref{analytic-transforms}.  The data, the uppermost time series in Fig.~\ref{analytic-transforms}a--c, is the eastward velocity recorded by a freely drifting subsurface oceanographic float \cite{richardson89-jpo,armi89-jpo,spall93-jmr}.  Such instruments are an important means of tracking the ocean circulation, and this and other such data may be downloaded from the World Ocean Circulation Experiment Subsurface Float Data Assembly Center (WFDAC) at \url{http://wfdac.whoi.edu}.  The oscillatory nature of the signal reflects the presence of an oceanic vortex [e.g, \citen{mcwilliams85-rvg}], properties of which may be inferred from the modulated oscillation as in \cite{lilly06-npg}.  This particular record has recently been used as an example in other studies \cite{rilling07-ispl,lilly10a-itsp}. More details regarding the data and its interpretation may be found in \cite{lilly10a-itsp}, but here we shall merely take it as a typical example of a modulated oscillation in noise.

\normalsize

Three wavelet transforms are shown in Fig.~\ref{analytic-transforms}d--f using the three wavelets shown in Fig.~\ref{analytic-morsies}a--c, together with locations of amplitude ridges. The resulting signal estimates and the differences between the original time series and the signal estimates are presented in Fig.~\ref{analytic-transforms}a--c.  The data is actually recorded as position, and the wavelet transform applied to this position record, but for clarity the time derivatives are presented in Fig.~\ref{analytic-transforms}a--c as these emphasize the oscillatory structure, rather than the lower-frequency meandering behavior which is also present.  The wavelet transform is taken at 74 logarithmically-spaced frequencies between radian frequencies $0.12$ and $2.39$; the data sample interval is one day.  All amplitude ridges are found whose length exceeds $2P_{\beta,\gamma}$---that is, $2\sqrt{4.5}=4.2$, $2\sqrt{9}=6$, and $2\sqrt{21}=9.2$ cycles for Fig.~\ref{analytic-transforms}a, b, and c respectively.  In all three cases there is only one such ridge, which extends nearly throughout the entire record.

All three of the estimates appear reasonable, and they are not drastically different from one another.  However, the residual curves in Fig.~\ref{analytic-transforms}a--c reveal some differences, with variability of the residual---and smoothness of the estimated signal---increasing as $P_{\beta,\gamma}$ increases.  Also, a major low-frequency fluctuation near time $t=0$ appears to have been missed by the smoothest estimate in Fig.~\ref{analytic-transforms}c.  An important issue is the ability to compare these signal estimates against one another to decide which is to be preferred; this is accomplished in Section~\ref{section:application} using the results of the subsequent development.  Discussion of the last two rows in Fig.~\ref{analytic-transforms} will be left until later.

\subsection{Outstanding Questions}
This section has presented essential elements of wavelet ridge analysis for modulated oscillations, as introduced by \cite{delprat92-itit} and extended by \cite{mallat}.  A number of questions may immediately be asked:
\begin{enumerate}
\item[(i)] What is the form of the time-varying bias terms in the ridge-based signal estimate?
\item[(ii)]  Are the amplitude and phase ridges the same, and if not which shows superior performance?
\item[(iii)] How should the wavelet properties be chosen in order to minimize bias?
\end{enumerate}
Addressing these questions is the goal of this paper, which we accomplish with special attention to the generalized Morse wavelets and using the data in  Fig.~\ref{analytic-transforms} as an example.

\section{Representation of Modulated Oscillations}\label{localsection}

In this section a local expansion of an analytic signal is constructed.  The hierarchy of terms in this expansion quantify increasingly higher-order local deviations of the signal from a constant amplitude, constant frequency sinusoid.

\subsection{A Local Representation}
\begin{definition}{Instantaneous Modulation Functions}\\
We aim to express the local variation of an analytic signal as a series of departures from a uniform oscillation. To this end we define functions $\widetilde\rho_n(t)$ for $n=1,2,\dots$ as
\begin{equation}
\widetilde\rho_n(t)\equiv \frac{1}{\omega^n(t)} \frac{1}{x_+(t)} \frac{d^n}{d\tau ^n}\left. \left[x_+(t+\tau)e^ {-i\omega(t)\tau}\right]\right|_{\tau =0}\label{demodulatedefinition}
\end{equation}
where $t$ on the right-hand-side is interpreted as a reference or ``global'' time, while $\tau$ is a ``local'' time.  The $\widetilde\rho_n(t)$ are the $\tau$-derivatives, evaluated at $\tau=0$, of $x_+(t+\tau)$ \emph{demodulated by} a uniform oscillation in local time $\tau$ having a frequency equal to the signal's instantaneous frequency $\omega(t)$ at the global time~$t$. Division by powers of $\omega(t)$ renders the $\widetilde\rho_n(t)$ dimensionless.  These functions will be called the \emph{instantaneous modulation functions}, and are to play an important part in what follows.  We are now in a position to state the following theorem.
\end{definition}
\begin{theorem}{The Local Modulation Expansion \label{Bell}}\\
Let $x_+(t)$ be an analytic signal defined by (\ref{analyticfilter}) for a real-valued $x(t)\in L^2({\mathbb{R}})$. Assume that $x_+(t)\in C^{N+1}\left[t,t+\tau\right]$, i.e. that $x_+(t)$ is $N+1$ times differentiable on the interval $\left[t,t+\tau\right]$ for some ``truncation level'' $N=0,1,2\dots$, and also that $x_+(t)\neq 0$ on that interval. Then $x_+(t+\tau)$ may then be expressed as an $N$th-order Taylor expansion in local time $\tau$
\begin{multline}
x_+(t+\tau)=   x_+(t)e^ {i\omega(t) \tau }\\
\left[1+\sum_{n= 1}^N \frac{1}{n!}\left[\omega(t)\tau\right]^n
\widetilde\rho_n( t)+R_{N+1}(\tau,t)\label{zBellN}
\right]
\end{multline}
where the form of the residual term
\begin{equation}
R_{N+1}(\tau,t)\equiv \frac{\left[\omega(t')\tau\right]^ {N+1}}{(N+1)!}\,\widetilde\rho_{N+1}( t')\quad
t'\in\left[t,t+\tau\right]
\label{RNdef}
\end{equation}
is found by employing the Lagrange form of the remainder in Taylor's theorem \cite[p880]{abramowitz}.  When $x(t)$ is an oscillatory signal, this expansion of a demodulated version of $x_+(t)$ can be expected to converge much more rapidly than a direct Taylor expansion of $x_+(t)$ itself.
\end{theorem}
\begin{proof}
The local modulation expansion is the Taylor series expansion of the complex-valued function
\begin{equation}
\frac{x_+(t+\tau)e^ {-i\omega(t) \tau }}{x_+(t)}\nonumber
\end{equation}
with respect to the variable $\tau$ about the point $\tau=0$.
\end{proof}

In the vicinity of some global time $t$, (\ref{zBellN}) represents the variation of $x_+(t+\tau)$ with respect to local time $\tau$ as a series of departures from a pure complex oscillation at the fixed frequency $\omega(t)$.  The $n$th-order instantaneous modulation function $\widetilde\rho_n(t)$ thus gives the contribution to the deviation of the signal from a pure sinusoid at $n$th order in the dimensionless local time $\omega(t)\tau$.  The instantaneous modulation functions are interpretable as fundamental quantities describing the deviation of the signal from a uniform oscillation.  For a constant-amplitude, constant-frequency sinusoid $x(t)=a_o\cos(\omega_o t)$, the instantaneous modulation functions of all orders are well-defined and vanish identically everywhere.  Signals which are considered ``oscillatory'' in the vicinity of some time $t$ should therefore have the instantaneous modulation functions $\widetilde\rho_n(t)$ not being too large in magnitude.

\subsection{The Instantaneous Modulation Functions}\label{section:modfun}

It remains to find the form of instantaneous modulation functions. In the following, we find it convenient to group the instantaneous frequency and bandwidth into a single complex-valued quantity
\begin{equation}
\eta(t)\equiv \omega(t)-i\upsilon(t)=-i\frac{d}{dt}\ln x_+(t)
\label{etadef}
\end{equation}
which we term the \emph{complex instantaneous frequency}.  The definition (\ref{etadef}) emphasizes that $\omega(t)$ and $\upsilon(t)$ are related as the imaginary and real parts, respectively, of the time derivative of $\ln x_+(t)$.  Since $\omega(t)$ and $\upsilon(t)$ often occur together as the complex-valued quantity~$\eta(t)$, this grouping will simplify subsequent expressions. We chose to multiply by $-i$ in \eqref{etadef} to make the definition concur with the usual instantaneous frequency for a pure sinusoidal function experiencing no amplitude modulation.

To find expressions for the instantaneous modulation functions, we first note that the bandwidth may be expressed as
\begin{equation}
\upsilon(t)=\frac{d}{d\tau}\left. \ln\left[x_+(t+\tau)e^ {-i\omega(t)\tau}\right]\right|_{\tau =0}    \end{equation}
while the complex instantaneous frequency $\eta(t)$ has an $n$th derivative given by
\begin{equation}
i\eta^{(n-1)}(t)=\frac{d^n}{d\tau ^n}\left. \ln\left[x_+(t+\tau)e^ {-i\omega(t)\tau}\right]\right|_{\tau =0},\quad n>1.
\end{equation}
Now assuming $\ln x_+(t)$ to be infinitely differentiable in the neighborhood of time~$t$,  (\ref{zBellN}) for infinite $N$ may be rewritten~as
\begin{multline}
\exp\left\{\displaystyle \ln \left[x_+(t+\tau)e^ {-i\omega(t) \tau }\right]-\ln x_+(t)\right\}
=\\1+\sum_{n= 1}^\infty \frac{1}{n!}\left[\omega(t)\tau\right]^n
\widetilde\rho_n( t)\label{zBellN2}
\end{multline}
which becomes, upon Taylor-expanding the exponent of the left-hand side,
\begin{multline}
\exp\left\{\upsilon(t)\tau+ \sum_{m=2}^\infty \frac{1}{m!}\,i\eta^{(m-1)}(t)\tau^m\right\}
=\\1+\sum_{n= 1}^\infty \frac{1}{n!}\left[\omega(t)\tau\right]^n
\widetilde\rho_n( t).\label{zBellN3}
\end{multline}
This expression indicates a relationship between the instantaneous modulation functions, appearing on the right-hand side, and the instantaneous bandwidth and derivatives of the complex instantaneous frequency on the left-hand side.

\subsection{Expressions Using Bell Polynomials}
To derive closed-form expressions for the instantaneous modulation functions, we turn to a special set of functions called the complete Bell polynomials \cite{bell33-aom,wiki:bell}. The complete Bell polynomial $B_n$, operating on $n$ arguments $c_1,c_2\,\ldots,c_n$, is defined to give the coefficients appearing in the expansion
\begin{equation}
\exp\left\{\sum_{n=1}^\infty \frac{1}{n!}\,c_n \tau^n\right\} =
\sum_{n=0}^\infty \frac{1}{n!}\,\tau^n B_n(c_1, c_2,\ldots,c_n )\label{zBelldefinition}
\end{equation}
with $B_0 \equiv 1$.  Expressions for the first four Bell polynomials
\begin{eqnarray}
\!\!\!\!\!\!\!\!\!B_1(c_1)& \seq&  c_1\label{bell1} \\
\!\!\!\!\!\!\!\!\!B_2(c_1,c_2)&  \seq&  c_1^ 2+c_2\label{bell2} \\
\!\!\!\!\!\!\!\!\!B_3(c_1,c_2,c_3)&  \seq&  c_1^3+3c_1c_2 +c_3\label{bell3}\\
\!\!\!\!\!\!\!\!\!B_4(c_1,c_2,c_3,c_4)& \seq&  c_1^4+6c_1^ 2c_2+4c_1c_3+3c_2^ 2+c_4\label{bell4}
\end{eqnarray}
can be verified directly by expanding (\ref{zBelldefinition}) and equating powers of  $\tau$ between the left-hand and right-hand sides.  More generally, the Bell polynomials satisfy a recursion relation \cite{bell33-aom}
\begin{equation}
B_{n}(c_1, c_2,\ldots,c_{n} ) =\sum_{p   = 0} ^ {n-1} \left(\!\!\begin{array}{c}n-1\\p \end{array}\!\!\right)c_{n -p} \,B_{p }(c_1, c_2,\ldots,c_{p} )\label{bellrecursion}
\end{equation}
for $n\geq 1$ given any $c_1, c_2,\ldots,c_{n}$.

Comparing (\ref{zBellN3}) with (\ref{zBelldefinition}) we find
 \begin{equation}
\widetilde\rho_n(t) =  B_n\left(\frac{\upsilon(t)}{\omega(t)}, \frac{i\eta'(t)}{\omega^ 2(t)},\ldots,\frac{i\eta^{(n-1)}(t)}{\omega^ {n}(t)}\right)\label{Belldefinition}
\end{equation}
as an expression for the $n$th instantaneous modulation function in terms of the $n$th-order Bell polynomial operating on the bandwidth $\upsilon(t)$ and the first $n-1$ derivatives of $\eta(t)$.  From (\ref{bell1}--\ref{bell3}) one then obtains
\begin{eqnarray}
\widetilde\rho_1(t)& =&  \frac{\upsilon(t)}{\omega(t)}\label{Bell1}\\
\widetilde\rho_2(t)& =&  \frac{\upsilon^ 2(t)}{\omega^ 2(t)}+\frac{i\eta'(t)}{\omega^ 2(t)}\label{Bell2} \\
\widetilde\rho_3(t)& = &  \frac{\upsilon^3(t)}{\omega^3(t)}+3\frac{\upsilon(t)}{\omega(t)}\frac{i\eta'(t)}{\omega^ 2(t)} +\frac{i\eta''(t)}{\omega^3(t)}\label{Bell3}
\end{eqnarray}
as the first three instantaneous modulation functions.

The $n$th instantaneous modulation function  $\widetilde\rho_n(t)$ thus combines powers of the bandwidth $\upsilon(t)$ and powers of time derivatives of the complex instantaneous frequency $\eta(t)$ into a measure of the $n$th order departure of the signal from a uniform oscillation.  On account of the nondimensionalization by powers of $\omega(t)$, we can interpret the rates of change of the amplitude and phase involved in $\widetilde\rho_n(t)$ to be on time scales proportional to the local instantaneous period $2\pi/\omega(t)$.  The first of these, $\widetilde\rho_1(t)$, is simply a nondimensional form of the bandwidth $\upsilon(t)$.  In general $\widetilde\rho_n(t)$ is complex-valued for~$n>1$.

\subsection{Examples}

As an example, consider the complex instantaneous frequency specified by
\begin{equation}
\eta(t)=\omega(t)-i\upsilon(t)=\omega_o\left[1-re^ {i\omega_1t}\right]
\end{equation}
where $\omega_o$ and  $\omega_1$ are real-valued constants, and $r$ is a potentially complex-valued constant with $|r|<1$. The normalized $n$th-order derivative of the complex instantaneous frequency~is
\begin{equation}
\left|\frac{\eta^ {(n)}(t)}{\omega^{n+1}(t)}\right|=\left|r\right|\frac{\left|\omega_o\right|\left|\omega_1\right|^n}{\left|\omega(t)\right|^{n+1}}
=\left|r\right|\frac{\left|\omega_o\right|}{\left|\omega(t)\right|}\left|\frac{\omega_1}{\omega(t)}\right|^{n}.
\end{equation}
This quantity decreases with increasing $n$ whenever $\left|\omega_1\right|$, the frequency at which the complex instantaneous frequency oscillates, is smaller than the local instantaneous frequency $\left|\omega(t)\right|$ itself.  If however $\left|\omega_1\right|$ exceeds $\left|\omega(t)\right|$, then derivatives of the complex instantaneous frequency will grow with $n$, eventually becoming non-negligible no matter how small one takes $\left|r\right|$.  Rapid fluctuations of the instantaneous frequency or bandwidth therefore cause the instantaneous modulation functions to fail to decay with increasing $n$.

As another example we return to the data analyzed in Fig.~\ref{analytic-transforms}. Panels (g--i) present the instantaneous frequency $\omega(t)$ and bandwidth $\upsilon(t)$ associated with the three estimated modulated oscillations in (a--c), while the corresponding values of $\widetilde\rho_2(t)$ and $\widetilde\rho_1^2(t)$ are shown in (j--l).  The instantaneous modulation functions $\widetilde\rho_1(t)$ and $\widetilde\rho_2(t)$ reveal the nature and magnitude of the first- and second-order modulation of the three estimated signals.  The most dramatic change is the increasing degree of smoothness, corresponding to the use of longer-duration wavelets, as one proceeds from the left column to the right column. The degree of amplitude variability increases in the latter half of all estimates, where the frequency has also increased, but we note from (g--i) that the amplitude modulation rate is generally very small compared to the instantaneous frequency, i.e.~$\upsilon(t)<<\omega(t)$.  Now, writing out (\ref{Bell2}) one finds
\begin{equation}
\widetilde\rho_2(t)
=\widetilde\rho_1^2(t)
+\frac{\upsilon'(t)}{\omega^ 2(t)}+\frac{i\omega'(t)}{\omega^ 2(t)}
\end{equation}
which shows that $\widetilde\rho_1^2(t)$ is a contributor to $\widetilde\rho_2(t)$.  Since panels (g--i) show that the bandwidth is very small compared to the instantaneous frequency, it is not surprising to find in (j--l) that $\widetilde\rho_2(t)$ contains only a minor contribution from $\widetilde\rho_1^2(t)=\upsilon^2(t)/\omega^2(t)$.  Instead we find that the real and imaginary parts of $\widetilde\rho_2(t)$ are roughly of an equal magnitude, implying comparable contributions from $\upsilon'(t)$ and $\omega'(t)$ in all three estimates.

\subsection{Signal Variability} \label{section:signalvar}

Using the instantaneous modulation functions we may now quantify the degree of departure of a signal from a uniform oscillation.

\begin{definition}{Local Signal Stability Level\label{stabilityassumption}}\\
Choose a truncation level $N_T$ for the local modulation expansion which is fixed  over some time interval $T$. Then the stability level $\delta_{N_T}$ of an analytic signal $x_+(t)$ is defined as the smallest positive constant which satisfies
\begin{eqnarray}
\left|\frac{\upsilon(t)}{\omega(t)}\right|&\le& \delta_{N_T}\quad\forall\quad t\in T \label{bandwidthcondition}\\
\left|\frac{\eta^{(n-1)}(t)}{\omega^n(t)}\right|&\le&\delta_{N_T}^n\quad\forall\quad t\in T, \quad 2\leq n\leq N_T.\label{derivativecondition1}
\end{eqnarray}
It is clear from (\ref{Bell1}--\ref{Bell3}), together with the recursive form of the Bell polynomials (\ref{bellrecursion}), that these conditions imply
\begin{eqnarray}
\widetilde\rho_n (t)&=&O(\delta_{N_T}^n)\quad\quad t\in T \quad 1 \leq n\leq N_T  \label{smallrho}
\end{eqnarray}
with powers of the bandwidth contributing at the same order as derivatives of the complex instantaneous frequency.  Furthermore, note that (\ref{bandwidthcondition}) and (\ref{derivativecondition1}) also imply
\begin{equation}
\frac{1}{\omega(t)}\frac{d}{dt}\,\widetilde\rho_n(t)=n \times O(\delta_{N_T}^{n+1})\quad\quad t\in T\quad  1 \leq n\leq N_T  \label{smalldrho}
\end{equation}
because, for example,
\begin{multline}
\frac{1}{\omega(t)}\frac{d}{dt}\,\widetilde\rho_2(t)=
\frac{1}{\omega(t)}\frac{d}{dt}\left[\frac{\upsilon^ 2(t)+i\eta'(t)}{\omega^ 2(t)}\right]\\
 = \frac{2\upsilon(t)\upsilon'(t)+i\eta''(t)}{\omega^ 3(t)}+O(\delta_{N_T} ^{4})= 2\times O(\delta_{N_T} ^{3})\label{rho2derivative}
\end{multline}
and similarly for higher-order $n$.
\end{definition}

The local stability level $\delta_{N_T}$ is determined by the variability of the signal, and may be different for different choices of truncation level $N_T$.  Thus the local stability level $\delta_{N_T}$ is a single number describing the extent to which \emph{any} square-integrable, $N_T+1$ times differentiable analytic signal $x_+(t)$ departs from a uniform oscillation at up to and including $N_T$th order.  When $\delta_{N_T}\ll 1$ it will be possible to obtain a greatly simplified representation of the AWT.  This is the key to obtaining direct closed-form expressions for the effect of signal modulation on the AWT and the ensuing ridge-based signal estimates, as we address in the next section.

\subsection{Wavelet Suitability Criteria}\label{section:suitabilty}

The local stability level can be used to match a wavelet to a signal in such a way that the wavelet ridge analysis will yield an accurate estimate of the signal.  The rational behind the ``wavelet suitability criteria'', which we now define, will become more clear when we apply them to the analytic wavelet transform in the next section.   These conditions constrain the choice of wavelet appropriate for a given oscillatory signal.

\begin{assumption}{Wavelet Suitability Criteria \label{waveletassumptions}}\\
Consider a signal characterized by local stability level $\delta_{N_T}$ with truncation level $N_T$ over some time interval $T$.  Given these quantifications of the signal's variability, we match a wavelet to the signal as follows.  We assume that the frequency-domain derivatives of the wavelet satisfy the following criteria at the peak frequency for $n\ge 2$
\begin{eqnarray}
\delta^{n/2}_{N_T}\frac{\widetilde\Psi_{n}(\omega_\psi )}{n!}& \le&
1  \quad \quad  \frac{n}{2}\in{\mathbb{N}}\label{waveletstability1}\\
\delta^{(n-1)/2}_{N_T}\frac{\widetilde\Psi_{n}(\omega_\psi )}{n!}& \le &
1  \quad \quad \frac{n-1}{2}\in{\mathbb{N}}.
\label{waveletstability2}
\end{eqnarray}
When $\delta_{N_T}$ is small, the quantities $(1/n!)\widetilde\Psi_{n}(\omega_\psi )$ are permitted to grow with increasing $n$ since then the powers of $\delta_{N_T}$ become smaller with increasing $n$.  The lowest-order suitability criteria, at $n=2$, implies that $P_\psi=|\widetilde\Psi_{2}|^{1/2}\le\sqrt {2/\delta_{N_T}}$.  This means that the time-domain wavelet is constrained so that is not too long---i.e. does not contain too many oscillations---compared to the degree of time variability of the signal, or alternatively, that the frequency-domain wavelet is constrained so that it is not too localized about the peak frequency.
\end{assumption}

Note that (\ref{waveletstability1}) and  (\ref{waveletstability2}) place a tighter condition on odd moments than on even moments.  This will simplify our analysis, and is reasonable because one expects the odd moments---which quantify the degree to which the wavelet departs from symmetry about its peak frequency---will be small for useful wavelet functions, as discussed earlier.   Since $\widetilde\Psi_{1}(\omega_\psi )$ vanishes by the definition of $\omega_\psi$, the lowest-order odd derivative to which these conditions apply is the third-order quantity $\widetilde\Psi_{3}(\omega_\psi )$.

\subsection{Suitability of the Generalized Morse Wavelets}\label{section:morsesuitabilty}

\begin{figure}[t]
\begin{center}
\includegraphics[width=3.5in,angle=0]{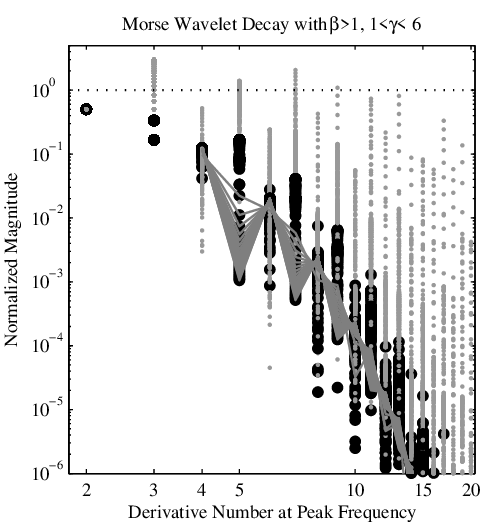}
\end{center}
            \caption{Decay of normalized frequency-domain derivatives of the generalized Morse wavelets, as discussed in the text. Each point shows the normalized $n$th frequency-domain derivative for a wavelet $\psi_{\beta,\gamma}(t)$  with integer $\beta$ in the range 1--21, and integer $\gamma$ in the range 1--11.  Values for $1\le\gamma\le 6$ are marked by large black circles, with values for $\gamma>6$ being gray dots.  The heavy gray lines show $\gamma=3$, beginning at $n=4$ since the $n=3$ terms vanish for these wavelets.
            }\label{derivative_decay}
\end{figure}

Next we find a range of parameter space for which the generalized Morse wavelets satisfy the wavelet suitability criteria.  Let us choose the wavelet such that $P_{\beta,\gamma}=\sqrt{ 2/\delta_{N_T}}$ for some local stability level $\delta_{N_T}$.  Thus for a given $\delta_{N_T}$, this fixes a curve in $(\beta,\gamma)$ space along which the lowest-order ($n=2$) suitability criterion is satisfied, and we must ask where along this curve the higher-order suitability criteria are also satisfied.  In Fig.~\ref{derivative_decay} we plot
\begin{eqnarray*}
\left(P_{\beta,\gamma}^2/2\right)^{-n/2}\frac{\widetilde\Psi_{n;_{\beta,\gamma}}(\omega_{_{\beta,\gamma}})}{n!} && \quad \frac{n}{2}\in{\mathbb{Z}}\label{newdecaycondition1}\\
\left(P_{\beta,\gamma}^2/2\right)^{-(n-1)/2}\frac{\widetilde\Psi_{n;_{\beta,\gamma}}(\omega_{_{\beta,\gamma}})}{n!}& &\quad \frac{n-1}{2}\in{\mathbb{Z}}\label{newdecaycondition2}
\end{eqnarray*}
for different values of the doublet $(\beta,\gamma)$, with $\beta\ge 1$ and $\gamma\ge 1$ and for $n\ge2$. If these two quantities are less than unity and with the choice $P_{\beta,\gamma}=\sqrt{ 2/\delta_{N_T}}$, then comparison with (\ref{waveletstability1}) and (\ref{waveletstability2}) shows that the wavelet suitability criteria are satisfied.

A difference in behavior is seen in Fig.~\ref{derivative_decay} for $1\le\gamma\le 6$, and other values of $\gamma$.  For  $1\le\gamma\le 6$ the normalized wavelet derivatives are always less than unity and decay rapidly with increasing~$n$.  For other values of $\gamma$ this is not the case, and we see that unity is occasionally exceeded at $n=3,5,7,9$. Also, outside the region  $1\le\gamma\le 6$ the rate at which the plotted terms decay with increasing $n$ is noticeably slower.  Thus if we are presented with a signal characterized by a local stability level $\delta_{N_T}$ over some time interval $T$ and for some truncation level $N_T$, we can choose a generalized Morse wavelet to satisfy the wavelet suitability criteria by setting $P_{\beta,\gamma}\le\sqrt{ 2/\delta_{N_T}}$ and choosing any $(\beta,\gamma)$ pair with  $\beta>1$ and $1\le\gamma\le 6$.  An application to the data in Fig.~\ref{analytic-transforms} will be given later.

\section{Analysis of Modulated Oscillations}\label{transformsection}

The goal of this section is the derivation of an expression for the analytic wavelet transform (AWT) of a potentially highly variable signal $x(t)$ which makes explicit the interaction between the analytic signal $x_+(t)$ and the wavelet.

\subsection{Additional wavelet properties}
Measures of the wavelet time-domain support and long-time decay will be needed.  The energy fraction function
\begin{equation}
\alpha_{\psi}(L)=\frac{\int_{-L}^{L} \left|\psi(t)\right|^2\,dt}{\int_{-\infty}^{\infty} \left|\psi(t)\right|^2\,dt}\label{alphadef}
\end{equation}
gives the ratio of the wavelet energy in a time window of half-width $L$ to the total energy.  The energy fraction is inverted by the time support function $L_\psi(\alpha)$
\begin{equation}
L_{\psi}(\alpha)\equiv\alpha_\psi^{-1}(\alpha)\label{Ldef}
\end{equation}
where the exponent ``$-1$'' in this context denotes the inverse function. $L_\psi(\alpha)$ associates a wavelet half-width with a given energy fraction, such that $\alpha$ is the fraction of the wavelet energy inside the time window $|t|<L_\psi(\alpha)$.  The long-time decay of the wavelet is specified by
\begin{equation}
\left|\psi(t)/\psi(0)\right|\sim |t|^{-r_\psi}\label{rdef}
\end{equation}
for some constant $r_\psi>0$, which takes on the value $r_{\beta,\gamma}=\beta+1$
for the generalized Morse wavelets \cite{lilly09-itsp}.

\subsection{Theorem}

Using the results of the preceding section, we may state the following theorem, which gives the exact form of the AWT of a potentially highly variable signal with a general analytic wavelet.

\begin{theorem}{The AWT Representation Theorem\label{transform}}\\
Fix  $(t,s)\in
{\mathbb{R}}\times {\mathbb{R}}_+$, choose a truncation level $N_T\in\mathbb{N}$ such that the wavelet decay satisfies $r_\psi\ge N_T+2$, and also an energy fraction $\alpha$ specifying a time support $L_\psi(\alpha)$.  Assume that \[\ln\left[x_+(t)\right]\in C^{N_T+1}\left[t-sL_\psi(\alpha),t+s L_\psi(\alpha)\right]\]
which implies that $|x_+(t)|\neq 0$ over the same interval.  The AWT of the real-valued signal $x(t)$ is then
\begin{multline}
W_\psi(t,s)=
\frac{1}{2}x_+(t)\Psi^*(s\omega(t))\times\\
\left[1+\sum_{n= 1}^N \frac{(-i)^n}{n!}\widetilde \Psi_n^*(s\omega(t)) \,
\widetilde\rho_n\left(t\right)
+\varepsilon_{\psi,N+1}(t,s)\right]
\label{zBelltransform}
\end{multline}
where $\varepsilon_{\psi,N+1}(t,s)$ is a transform residual given by (\ref{newvarepsilondef}) in Appendix~\ref{appendix:representation}.
\end{theorem}
\begin{proof}
The proof is provided in Appendix~\ref{appendix:representation}, together with bounds on the transform residual.  For an intuitive illustration of the basic idea, here we prove an idealized special case.  In this paragraph we take $\psi(t)$ to be a filter that is exponentially decaying in time, which means it cannot be an analytic function; thus a term arising from non-analyticity of $\psi(t)$ emerges here but not in (\ref{zBelltransform}). The analytic signal is assumed everywhere infinitely differentiable and non-vanishing.  After a change of variables, and noting $x(t)=[x_+(t)+x_+^*(t)]/2$, the wavelet transform (\ref{wavetrans}) becomes
\begin{multline}
W_\psi(t,s)
=\frac{1}{2}\frac{1}{s}\int_{-\infty}^{\infty} \psi^*\left( \frac{\tau}{s}\right)\left[ x_ +(t+\tau)+x_ +^*(t+\tau)\right]\,d\tau
\\\equiv W_{\psi,x_+}(t,s)+W_{\psi,x_+^*}(t,s)
\label{wavetrans3}
\end{multline}
in which we have implicitly defined an analytic portion $W_{\psi,x_+}(t,s)$ and an anti-analytic portion $W_{\psi,x_+^*}(t,s)$.  For the former, substituting the local modulation expansion (\ref{zBellN}) gives
\begin{multline}
W_{\psi,x_+}(t,s) =\frac{1}{2}\,x_+(t)\frac{1}{s}\int_{-\infty}^{\infty} \psi^*\left( \frac{\tau}{s}\right) \,e^ {i\omega(t)  \tau}\times \\ \left[1+\sum_{n= 1}^\infty \frac{1}{n!}\left[\omega(t)\tau\right]^n\widetilde\rho_n(t)\right]\,d\tau.\label{WSigmainfinite}
\end{multline}
After the change of variables $\tau/s=u$, and exchanging the order of summation and integration, we have
\begin{multline}
  W_{\psi,x_+}(t,s)  = \frac{1}{2}\,x_+(t)\times \left[\Psi^*(s\omega(t))\right.\\\left.+\sum_{n= 1}^\infty \frac{(-i)^n}{n!} \widetilde\rho_n(t)\int_{-\infty}^{\infty}\left[is\omega(t)\tau\right]^n \psi^*(\tau) \,e^ {is\omega(t)  \tau} \,d\tau\right]\label{WSigmainfinite2}
\end{multline}
where we have introduced canceling factors of $i^n$ and $(-i)^n$.  Now, the $n$th dimensionless derivative (\ref{normalizedwavelet}) may be written as
\begin{equation}
\widetilde\Psi_n(\omega) =\frac{1}{\Psi(\omega) } \int_{-\infty}^{\infty} \label{nthderivative} \left(i\omega\tau\right)^n\psi(\tau) \,e^ {-i\omega \tau} \,d\tau
\end{equation}
by differentiating the Fourier representation of $\Psi(\omega)$. Substituting this into (\ref{WSigmainfinite2}) obtains the right-hand side of (\ref{zBelltransform}) with infinite~$N$.  Thus essentially (\ref{zBelltransform}) arises by noting that taking AWT of the Taylor-expanded, demodulated signal involves forming the time domain moments, hence the frequency domain derivatives, of the analyzing wavelet.  The proof in Appendix~\ref{appendix:representation} handles the truncation to a finite number of terms $N$ in the summation, together with complications arising from the polynomial rate of time decay of the wavelet.
\end{proof}

\subsection{Comments and Interpretation}

At this stage we offer some comments on the importance and interpretation of the preceding theorem, which we view as a fundamental result.  The AWT representation theorem (\ref{zBelltransform}) shows that the AWT is generated by the interaction of the frequency-domain derivatives of the wavelet with certain time-varying signal quantities---the instantaneous modulation functions introduced in the previous section. Each higher-order frequency-domain derivative of the wavelet $\widetilde\Psi_n(\omega)$ interacts with a higher-order measure $\widetilde\rho_n(t)$ of the variability of the signal.  The roles of amplitude and frequency modulation in setting transform properties are explicitly included. An advantage of the AWT representation theorem is that it permits us to compare different analytic wavelets for a given signal by comparing their frequency-domain derivatives.

In contrast to previous works \cite{delprat92-itit,mallat}, which assume that the signal bandwidth is small and that the bandwidth and instantaneous frequency are essentially constant, (\ref{zBelltransform}) resolves the hierarchy of nonlinear terms and is therefore useful for a much broader variety of local signal behavior. It can be seen as a substantial generalization of the pioneering work of Delprat et al. ~\cite{delprat92-itit} and Mallat~\cite{mallat}. In particular, Theorem~4.5 of Mallat~\cite{mallat} is roughly equivalent to (\ref{zBelltransform}) for the case $N=0$, that is, with the error term including everything except for the leading term of unity.  Mallat's derivation assumes a particular form for the wavelet---a real-valued envelope multiplied by a complex exponential---which cannot be strictly analytic, but which mimics the form of the popular Morlet wavelet \cite{holschneider}. The original proof by Delprat et al. \cite{delprat92-itit} of a result related to Mallat's relied on a stationary phase approximation, and similarly required the assumption of negligible modulation for both the wavelet and the signal.

\subsection{Compression Along Instantaneous Frequency Curves}\label{compressionsection}

Here we use the AWT representation theorem to examine the wavelet transform along the instantaneous frequency curve, a key theoretical quantity which controls the behavior of the ridge-based signal estimator. This sets the stage for the application to wavelet ridge analysis in the next section.

\begin{definition}{The Localized Analytic Signal}\\
Evaluating the AWT along the instantaneous frequency curve yields a fundamental object reflecting the joint properties of the signal and the wavelet,
\begin{equation}
 x_\psi  (t) \equiv  W_\psi\left(t, \omega_\psi/\omega(t)\right) \label{locallyanalyticdefined}
\end{equation}
which we term the \emph{localized analytic signal}. From the AWT representation theorem (\ref{zBelltransform}), we find immediately
\begin{multline}
 x_\psi  (t) =x_+(t)\times\\\left[1+\sum_{n= 2}^N  \frac{(-i)^n}{n!} \widetilde\Psi^*_{n}(\omega_\psi)\widetilde\rho_n(t)\label{locallyanalytic}
+\epsilon_{\psi,N+1}\left(t\right)\right]
\end{multline}
[recalling $\Psi(\omega_\psi)\equiv 2$] where the residual in this expression
\begin{eqnarray}
\epsilon_{\psi,N+1}\left(t\right) &\equiv &\varepsilon_{\psi,N+1}\left(t,\omega_\psi/\omega(t)\right)
\end{eqnarray}
is the transform residual appearing in (\ref{zBelltransform}) evaluated along the instantaneous frequency curve.  Note that no $n=1$ term appears in (\ref{locallyanalytic}) due to the fact that $\widetilde\Psi^*_{1}(\omega_\psi)=0$ by definition.
\end{definition}

The localized analytic signal $x_\psi(t)$ is a non-uniform and nonlinear filtering of the signal by the wavelet. It can be seen as a sequence of local projections of the signal onto a set of analyzing functions, in which the analyzing functions---the rescaled wavelets $\psi(t/s)/s$---are scaled to be proportional to the local instantaneous period.  The localized analytic signal represents a non-uniform filtering because the analyzing wavelet changes in scale across time, and this filtering is nonlinear since the scale of the wavelet depends upon the instantaneous frequency of the signal being analyzed.  It is clear that the localized analytic signal reduces to a linear time-invariant filtering when the instantaneous frequency is constant, since then it can be considered as merely the result of a convolution of the signal with some fixed wavelet function.

In general, the localized analytic signal is not itself precisely analytic. Its analyticity is compromised on account of the localization. Taking the Fourier transform of  (\ref{locallyanalytic}), and noting that time-domain multiplications become frequency-domain convolutions, we find that the Fourier transform of $x_\psi (t)$ is given by
\begin{multline}
X_\psi(\omega)=X_+(\omega)+\sum_{n=2}^{N}\frac{(-i)^n
\widetilde{\Psi}_n^*(\omega_\psi )}{n!}\times \\\frac{1}{2\pi}\int_{-\infty}^{\infty}X_+(\omega')
\widetilde P_n(\omega-\omega')\,d\omega' +\calE_{\psi,N+1}(\omega)
\end{multline}
where $\widetilde P_n(\omega)$ and $\calE_{\psi,N+1}(\omega)$ are defined as the Fourier transforms of $\widetilde\rho_n(t)$ and of the product $x_+(t)\epsilon_{\psi,N+1}(t)$, respectively. If $\widetilde\rho_n(t)=
0$ for all $n\ge 2$, then $X_\psi(\omega)$ has no support on negative frequencies; otherwise the convolutions may distribute energy to negative frequencies.  In this manner the interaction between signal and wavelet can cause the localized analytic signal $x_\psi (t)$ to deviate from exact analyticity.

\subsection{Bias of the Localized Analytic Signal}\label{section:localizedbias}

Making use of the local stability level, we may match the wavelet to the signal in such a way that the difference between  $x_\psi(t)$ and $x_+(t)$ is small.  The localized analytic signal (\ref{locallyanalytic}) involves a hierarchy of interactions between the instantaneous modulation functions and the frequency-domain derivatives of the wavelet.  In order that the localized analytic signal $x_\psi(t)$ be close to the true analytic signal $x_+(t)$, these interaction terms must be kept small.  At this point we invoke the wavelet suitability criteria (\ref{waveletstability1}) and (\ref{waveletstability2}), introduced in Section~\ref{section:suitabilty}, to obtain a simple expression for the bias of the localized analytic signal.

Assume the signal is characterized by local stability level $\delta_{N_T}$ for a truncation level $N_T$ over a time interval $T$.  Invoking the wavelet suitability criteria, the difference between the localized analytic signal and true analytic signal is for $t\in T$
\begin{multline}
\Delta x_\psi  (t) \equiv \frac{x_\psi(t) -x_+(t) }{x_+(t)}=
\overset{O(\delta_{N_T})}{\overbrace{-\frac{1}{2}\widetilde\Psi^*_{2}(\omega_\psi)\widetilde\rho_2(t)}} \\\overset{O(\delta_{N_T}^2)}{\overbrace{-\frac{i}{3!}\widetilde\Psi^*_{3}(\omega_\psi)\widetilde\rho_3(t)+\frac{1}{4!}\widetilde\Psi^*_{4}(\omega_\psi)\widetilde\rho_4(t)}}
\label{locallyanalyticerror}
+\epsilon_{\psi,5}\left(t\right)
\end{multline}
[using (\ref{zBelltransform})] where we have set the truncation level set to $N_T=4$.  Thus with a suitable choice of wavelet, the deviation of the localized analytic signal consists of a series of terms representing increasingly higher-order interactions of the signal with the wavelet, which diminish with increasing order.

Ensuring that a small value of $\Delta x_\psi  (t)$ is obtained places two constraints the choice of wavelet.  Firstly, as discussed in Section~\ref{section:suitabilty}, $-\frac{1}{2}\widetilde\Psi^*_{2}(\omega_\psi)=P_\psi^2$ is a measure of the (squared) wavelet duration and should be chosen so that the magnitude of the product $\widetilde\Psi^*_{2}(\omega_\psi)\widetilde\rho_2(t)$ is small.  Note that this lowest-order contribution to $\Delta x_\psi (t)$, at first order in $\delta_{N_T}$, is associated with $\widetilde\rho_2(t)$ rather than with $\widetilde\rho_1(t)$; this arises on account of the vanishing of $\widetilde\Psi_{1}(\omega)$ at the peak frequency $\omega_\psi$.  Secondly, it is important to make an appropriate choice of wavelet with fixed $P_\psi$ so that the higher-order terms are small.  For example, Fig.~\ref{analytic-morsies} shows that large absolute values of $\widetilde\Psi_{3}(\omega_\psi)$ correspond to a high degree of frequency-domain asymmetry; thus the $n=3$ suitability criterion can be interpreted as a bound on an acceptable degree of asymmetry.  More generally, if the suitability criteria are satisfied, contributions from higher-order instantaneous modulation functions appear at higher orders in $\delta_{N_T}$ than the leading term involving the duration $P_\psi$, and may consequently be neglected when $\delta_{N_T}$ is sufficiently small.

It is instructive to examine the form of the amplitude and phase of the localized analytic signal if we keep only the lowest-order term in the expansion.  With $N_T=2$ we have
\begin{equation}
x_\psi(t) =x_+(t)\left[1+\frac{1}{2}P_\psi^2\widetilde\rho_2(t) + \epsilon_{\psi,3}(t)\right]\label{simplifiedlocalized}
\end{equation}
recalling the definition (\ref{pdef}) of $P_\psi$. The amplitude and phase of the localized analytic signal are implicitly defined via
\begin{equation}
a_\psi(t)e^{i \phi_\psi(t)}\equiv x_\psi(t) \label{localizedexpansion}
\end{equation}
which, it should be pointed out, are not necessarily a canonical pair because $x_\psi(t)$ is not necessarily precisely analytic.  We may introduce the deviations
\begin{eqnarray}
 a_\psi(t) &=& a_+(t) \left[1+\Delta a_\psi(t)\right]\\
 \phi_\psi(t) &=& \phi_+(t)+\Delta \phi_\psi(t)
\end{eqnarray}
and inserting these expressions into (\ref{localizedexpansion}), one obtains
\begin{multline}
x_\psi(t)
=x_+(t) \left[1+\Delta a_\psi(t) +i\Delta \phi_\psi(t) +\right.\\\left. O\left\{\Delta \phi_\psi^2(t)\right\}+O\left\{\Delta a_\psi(t)\Delta \phi_\psi(t)\right\}
\right].
\end{multline}
Equating terms with (\ref{simplifiedlocalized}) then leads to
\begin{eqnarray} \!\!\!\!\!\!\!\!\!\!\!\!  \Delta a_\psi(t)& = &
\frac{1}{2}\label{locallyamplitude} \frac{P_\psi^2}{\omega^2(t)}\,\frac{a_+''(t)}{a_+(t)}
+O\left\{\epsilon_{\psi,3}(t)\right\}
+O\left(\delta_{N_T}^2\right)\\
\!\!\!\!\!\!\!\!\!\!\!\!  \Delta \phi_\psi(t) & = &\label{locallyphase}
 \frac{1}{2}\frac{P_\psi^2}{\omega^2(t)}\,\phi_+''(t)+
 O\left\{\epsilon_{\psi,3}(t)\right\}+O\left(\delta_{N_T}^2\right)
\end{eqnarray}
for the amplitude and phase deviation, respectively.

In both cases the deviation is proportional to the second derivative, or \emph{curvature}, of the quantity of interest.  Since $P_\psi$ is a measure of the time duration of the wavelet, these deviations are small when the amplitude curvature and phase curvature, respectively, are small over the time support of the wavelet.  Both amplitude and phase are underestimated during local maxima and overestimated during local minima.  This is intuitive behavior since the localized analytic signal is essentially a smoothed version of the true analytic signal.

In summary, to minimize the bias of the amplitude and phase of the localized analytic signal, one should make the magnitude of $P_\psi^2 \widetilde\rho_2(t)$ small, while simultaneously enforcing the wavelet suitability criteria.   This raises the question of what is the smallest useful value for $P_\psi$, an issue which will be addressed Section~\ref{section:application}.

\section{Estimation of Modulated Oscillations}\label{ridgesection}

Building on the results of the previous sections,  we explicitly identify the hierarchy of time-dependent bias terms associated with estimation of analytic signal properties using the wavelet ridge method, and show how to choose wavelets which minimize these terms.

\subsection{Transform Near an Instantaneous Frequency Curve}\label{section:transformnear}
In the preceding section we defined the localized analytic signal, which is simply the set of values taken by the wavelet along the instantaneous frequency curve.  In practice, however, this quantity is not known because the instantaneous frequency curve is not known.  However the ridge curves defined in Section~\ref{section:ridgeintro} may be identified directly from the transform, and these will be shown to closely approximate an instantaneous frequency curve. The ridge-based signal estimate is then found by evaluating the wavelet transform along a ridge:
\begin{equation*}
\widehat x^{\{\cdot\}}_{+,\psi}(t)\equiv W_\psi\left( t, s^{\{\cdot\}}(t)\right)\quad t\in  T^{\,\{\cdot\}} \quad \quad   (\ref{estimator})
\end{equation*}
where the ``$\cdot$'' is either an ``$a$'' or a ``$p$'' to refer to an amplitude ridge or a phase ridge, respectively, and where $T^{\,\{\cdot\}}$ is the time interval over which the ridge exists in practice.

In order to obtain an expression for the bias of the ridge-based signal estimate, we must account for the deviation of the ridge $s^{\{\cdot\}}$ from the instantaneous frequency curve $\omega_\psi/\omega(t)$.  To this end we employ another Taylor-series expansion and express the wavelet transform along the ridge in terms of the wavelet transform along the instantaneous frequency curve.  There emerge powers of a quantity
\begin{equation}
\Delta\omega(t,s)\equiv\frac{s\omega(t)}{\omega_\psi }-1
\end{equation}
which we refer to as the \emph{scale derivation} since it gives the departure of a given scale from the instantaneous frequency curve.  We then obtain the following result.

\begin{theorem}{The AWT Scale Deviation Expansion}\\
With the same assumptions as the AWT representation theorem, and the additional assumption that $\Psi(\omega)\in C^{\infty}\left(0,\infty\right)$, the AWT representation theorem (\ref{zBelltransform}) can be cast in the form
\begin{multline}
W_\psi(t,s)= x_+(t)\left[\sum_{m=0}^{\infty} \sum_{n=0}^{N} \sum_{p=0}^n \frac{(-i)^n}{(n-p)!p!
m!}\times \right.\\\left.
\widetilde{\Psi}_{m+n}^*(\omega_\psi )\widetilde\rho_n(t)
\left[\Delta\omega(t,s)\right]^{m+p}+\varepsilon_{\psi,N+1}(t,s)\right]\label{triplesummation}
\end{multline}
where now we define $\widetilde\rho_0(t)\equiv 1$ for convenience.  This expansion involves a triple summation over orders of the wavelet derivatives evaluated at the fixed frequency $\omega_\psi$, orders of the instantaneous modulation functions, and powers of the scale deviation.
\end{theorem}
\begin{proof}
The proof is given in Appendix~~\ref{section:deviation}.
\end{proof}

The AWT scale deviation expansion relates the AWT of the signal along the instantaneous frequency curve to the AWT at all scales. It can be seen as expressing how the compression of the signal along the instantaneous frequency curve extends across the time/scale plane.  Using the local stability level, and the wavelet suitability criteria, we can simplify the scale deviation expansion (\ref{triplesummation}) in the vicinity of an instantaneous frequency curve.  One final definition is also necessary.

\begin{definition}{The Instantaneous Frequency Neighborhood\label{neighborhood}}\\
The {\em instantaneous frequency  $n$-neighborhood} is defined by
\begin{equation}
{\cal R}_{n;\psi}(\delta_{N_T}  )=\left\{(t,s):\,\Delta\omega(t,s)
=O(\delta^n_{N_T}  )\right\}.\label{neighborhooddef}
\end{equation}
Thus, the $n$-neighborhood is a region of the time/scale plane within which the scale deviation is of $n$th order with respect to the local stability level $\delta_{N_T}$.  As $n$ increases, a region of the time/scale plane more tightly localized around $\omega_\psi/\omega(t)$ is specified.  This permits us to quantify the magnitude of the departure of a given scale point from the instantaneous frequency curve.  We can now state the following theorem.
\end{definition}

\begin{theorem}{The AWT Ridge Representation Theorem}\\\label{ridgetheoremtheorem}
Let $x_+(t)\in $ with $t\in T$ be characterized by stability level $\delta_{N_T}$ up to order $N_T$, and assume the wavelet is chosen such that the suitability criteria hold. Finally constrain the scale $s$ such that $(t,s)\in {\cal R}_{2;\psi}(\delta_{N_T})$, implying the scale deviation is of \emph{second} order in $\delta_{N_T}$.  We write the AWT of $x(t)$ as
\begin{equation}
\label{transformpartition}
 W_\psi(t,s) = x_+(t)\left[1+\Delta x_\psi(t)+\Delta W_\psi(t,s)\right]
\end{equation}
which serves to \emph{define} $\Delta W_\psi(t,s)$. Thus $ W_\psi(t,s)$ is separated into a scale-independent perturbation $\Delta x_\psi(t)$ and a scale-dependent perturbation $\Delta W_\psi(t,s)$.  The form of $\Delta x_\psi(t)$ was given earlier in (\ref{locallyanalyticerror}), while by subtraction from (\ref{triplesummation}) we find
\begin{multline}
\Delta W_\psi(t,s)=\overset{O(\delta_{N_T}^2)}{\overbrace{-i\Delta\omega(t,s)\widetilde\rho_1(t) \widetilde\Psi^*_{2}(\omega_\psi)}}\\
\overset{O(\delta_{N_T}^3)}{\overbrace{-\Delta\omega(t,s)\left[\widetilde\rho_2(t)\left(\widetilde\Psi^*_{2}(\omega_\psi)+\frac{1}{2}\widetilde\Psi^*_{3}(\omega_\psi)\right) -\frac{i}{3!}\widetilde\rho_3(t)\widetilde\Psi^*_{4}(\omega_\psi)\right]}}\\
\overset{O(\delta_{N_T}^3)}{\overbrace{+\frac{1}{2}\left[\Delta\omega(t,s)\right]^2\widetilde\Psi^*_{2}(\omega_\psi)}}+O\left(\delta^{4}_{N_T}\right)+ \varepsilon_{\psi,4}(t,s)\label{ridgetheorem}
\end{multline}
for the form of the scale-dependent perturbation $\Delta W_\psi(t,s)$.
\end{theorem}
\begin{proof}
This theorem follows directly from the AWT scale deviation expansion (\ref{triplesummation}) truncated at  $N= 3$,  together with the stated assumptions.
\end{proof}

The AWT ridge representation theorem gives the form of the analytic wavelet transform in the vicinity of an instantaneous frequency curve, explicitly resolving the effects of modulation up to third order in $\delta_{N_T}$.  An important point is that the scale-independent perturbation $\Delta x_\psi (t)$ contains the lowest-order term, at first order in $\delta_{N_T}$, while the scale-dependent perturbation $\Delta W_\psi(t,s)$ contains only terms of second order or higher in $\delta_{N_T}$.

Note that the signal stability level $\delta_{N_T}$ has been used in several different ways.  Its value is set from (\ref{bandwidthcondition}) and (\ref{derivativecondition1}) by the values of the derivatives of the analytic signal over some time interval $T$.  These conditions involve up to the $N_T$th derivative of the signal, or the $(N_T-1)$th derivative of the complex instantaneous frequency; here $N_T$ is a number that can be chosen, up to the degree of the signal's differentiability, and its choice will impact the value of $\delta_{N_T}$ that we find from the signal.  The signal stability level $\delta_{N_T}$ then constrains the choice of wavelet via the wavelet suitability criteria (\ref{waveletstability1}) and (\ref{waveletstability2}).  Finally, the local stability level is also involved in the notion of the instantaneous frequency neighborhood (\ref{neighborhooddef}), which indicates a region of the time-scale plane within which a simplification of the AWT representation theorem may be found.  The use of $\delta_{N_T}$ for the instantaneous frequency neighborhood has enabled an ordering of terms both on and off the instantaneous frequency curve using a single small parameter.

\subsection{Expressions for the Ridge-Based Signal Estimates}
Using the AWT ridge representation theorem, we can now obtain closed-form expressions for the ridge curves and the associated estimate of the analytic signal.   Henceforth, for simplicity, we assume that  $\widetilde\Psi_{n}(\omega_\psi)$ is real-valued for  $n\leq 4$, as is the case for the generalized Morse wavelet family of analytic wavelets \cite{olhede02-itsp,lilly09-itsp}.

In Appendix~\ref{section:ridgecurves} we find that both the amplitude ridge condition (\ref{ampridge}) and the phase ridge condition (\ref{phaseridge}) have unique solutions within the 2-neighborhood of the instantaneous frequency curve.  The amplitude ridges are found to have the explicit form
\begin{multline}
\widehat s^{\,\{a\}}(t)
=\frac{\omega_\psi}{\omega(t)}\label{amplitudecurve}
\left[1+ \left(\frac{\upsilon^2(t)}{\omega^2(t)} +\frac{\upsilon'(t)}{\omega^2(t)}\right) \left(1+\frac{1}{2}\frac{\widetilde\Psi_{3}(\omega_\psi)}{\widetilde\Psi_{2}(\omega_\psi)}\right)
\right.
\\\left.+\frac{1}{6}
\left(\frac{\omega''(t)}{\omega^3(t)}+3\frac{\upsilon(t)}{\omega(t)}\frac{\omega'(t)}{\omega^2(t)}  \right)\frac{\widetilde\Psi_{4}(\omega_\psi)}{\widetilde\Psi_{2}(\omega_\psi)}-\frac{1}{2}\frac{\upsilon(t)}{\omega(t)} \frac{\omega'(t)}{\omega^2(t)}  \widetilde\Psi_{2}(\omega_\psi)\right.
\\\left.+
O\left(\delta_{N_T} ^{3}\right)+O\left\{ \epsilon_{\psi,3}^{\{a,s\}}(t)\right\}\right]
\end{multline}
while the phase ridges are given by
\begin{multline}
\widehat s^{\,\{p\}}(t)
=\frac{\omega_\psi}{\omega(t)}\left[1
 +\left(\frac{1}{2} \frac{\omega''(t)}{\omega^3(t)}+ \frac{\upsilon(t)}{\omega(t)} \frac{\omega'(t)}{\omega^2(t)}  \right) \widetilde\Psi_{2}(\omega_\psi)\right.\\\left. +O (\delta_{N_T} ^{3})+O\left\{ \epsilon_{\psi,3}^{\{p,t\}}(t)\right\}\right]\label{phasecurve}
\end{multline}
where we have resolved terms up to second order in $\delta_{N_T}$.
Note truncation levels $N_T=4$ and $N_T=3$ have been used in the former and latter cases respectively.  The forms of the residual terms $\epsilon_{\psi,3}^{\{\cdot,\cdot\}}(t)$ are given by (\ref{epsilonsdef}) and (\ref{epsilontdef})  of Appendix~\ref{section:ridgecurves}.

These expressions have a rather complicated form, but there are two simple messages.   Firstly, amplitude and phase ridges are definitely \emph{not} the same, except perhaps for some special choices of wavelet properties. Secondly, from the wavelet suitability criteria it follows that all the resolved perturbation terms in (\ref{amplitudecurve}) and (\ref{phasecurve}) are of \emph{second} order in $\delta_{N_T}$; there are no terms at first order.   Thus both types of ridge curves are of the form
\begin{equation}
\widehat s^{\{\cdot\}}(t)\label{ridgesecondorder}
=\frac{\omega_\psi}{\omega(t)}\left[1+O (\delta_{N_T} ^{2}) +\ldots\right]
\end{equation}
where the ellipses indicate the omitted residual terms; here the superscript ``$\cdot$'' indicates either ``$a$'' or ``$p$''. This is important because the deviations of the ridge curves from the instantaneous frequency curve, and form each other, will be a higher-order effect compared to the first-order deviation of the localized analytic signal from the true analytic signal.

An expression for the ridge-based signal estimate is found by substituting (\ref{ridgesecondorder}) into (\ref{ridgetheorem}) for $\Delta W_\psi(t,s)$.  With a truncation level of $N_T=2$, one finds
\begin{multline}
\widehat x^{\{\cdot\}}_{+,\psi}(t)= \\x_+(t)\left[1+\frac{1}{2}P_{\psi}^2\widetilde\rho_2(t)+O \left(\delta_{N_T} ^{2}\right)+O\left\{ \epsilon_{\psi,3}^{\{\cdot,\cdot\}}(t)\right\}\right]\label{leadingerrorexpression}
\end{multline}
which we note is identical, apart from the final residual term, to expression (\ref{simplifiedlocalized}) for the localized analytic signal $x_\psi(t)$.  Thus the leading-order error term is due to the scale-independent departure of the localized analytic signal---i.e. the AWT along the instantaneous frequency curve---from the true analytic signal, rather than the departure of the estimated instantaneous frequency curve from the true instantaneous frequency curve.  We may alternately write (\ref{leadingerrorexpression}) as
\begin{equation}
\widehat x^{\{\cdot\}}_{+,\psi}(t)= x_\psi(t)\left[1+O \left(\delta_{N_T}^2\right)+O\left\{ \epsilon_{\psi,3}^{\{\cdot,\cdot\}}(t)\right\}\right]
\end{equation}
which states that the ridge-based signal estimate accurately recovers the \emph{localized analytic signal}.  The  difference between the two types of ridges will be of negligible importance when $\delta_{N_T}$ is small.

To this theoretical result, we should add a caveat.  While the perturbation analysis suggests there is no reason to prefer amplitude versus phase ridges, in practice we find the amplitude ridges to be superior.  When both exist, we generally find they are indeed very close to one another, as expected by the perturbation analysis, but the phase ridges have a greater tendency to ``break'' at isolated points where modulation is particularly strong.  In fact we find this to be the case when applying the phase ridge algorithm to the example in Fig.~\ref{analytic-transforms}  with identical settings as for the amplitude ridges (not shown).  Therefore based on experience we favor the amplitude ridges.  It is conceivable, however, that this difference in  performance is due to some aspect of our particular numerical implementation.

We may similarly find the amplitude and phase estimates associated with the ridge-based signal estimate.  Writing the estimated analytic signal in terms of an amplitude and phase
\begin{equation}
a^{\{\cdot\}}_{+,\psi}(t)e^{i \phi^{\{\cdot\}}_{+,\psi}(t)}\equiv x^{\{\cdot\}}_{+,\psi}(t)\label{localizedexpansionestimated}
\end{equation}
we find, following the development in Section~\ref{section:localizedbias}
\begin{eqnarray}
\widehat a^{\{\cdot\}}_{+,\psi}(t)&=& a_\psi(t)\left[1+O \left(\delta_{N_T} ^{2}\right)+O\left\{ \epsilon_{\psi,3}^{\{\cdot,\cdot\}}(t)\right\}\right]\\
\widehat \phi^{\{\cdot\}}_{+,\psi}(t)&=& \phi_\psi(t)\left[1+O \left(\delta_{N_T} ^{2}\right)+O\left\{ \epsilon_{\psi,3}^{\{\cdot,\cdot\}}(t)\right\}\right]
\end{eqnarray}
so that the estimated amplitude and phase are the same as the amplitude and phase of the localized analytic signal up to second order in $\delta_{N_T}$.

At this point we return to the assumption made in the derivation of the AWT Ridge Representation Theorem (\ref{transformpartition}) that the ridge points lie within a 2-neighborhood of an instantaneous frequency curve.  The solution to the ridge equations then gives (\ref{ridgesecondorder}) which is consistent with that assumption.  It turns out that, had we derived the AWT ridge representation theorem with the less restrictive assumption that scale $s$ lies within the 1-neighborhood of an instantaneous frequency curve, we would have again found (\ref{ridgesecondorder}) stating that the ridge curve in fact lies within the smaller 2-neighborhood.  This is why we have chosen to assume that the ridge points within the 2-neighborhood from the outset.

\subsection{Instantaneous Frequency and Bandwidth Estimation}

Expressions for the estimated instantaneous frequency and bandwidth can also be found.  A direct method of estimating the instantaneous frequency is simply through the scale frequency associated with the ridge curves, i.e.
\begin{equation}
\widehat \omega^{\{\cdot, s\}}(t)\equiv \frac{\omega_\psi}{\widehat s^{\{\cdot\}}(t)}=\omega(t)\left[1+O(\delta_{N_T}^2)\right]. \label{directestimate}
\end{equation}
However, instantaneous frequency estimates formed in this way are not very satisfactory because they reflect the discrete scale levels $s$ used in the numerical evaluation of the wavelet transform.  Likewise, one could differentiate the amplitude and phase of the estimated analytic signal $\widehat x_+(t)$, but in our experience the discrete implementation of this differentiation tends to lead to rather noisy estimates.

A better way to estimate both instantaneous frequency and bandwidth is to first form the quantities
\begin{eqnarray}
\Omega_\psi(t,s)\equiv   \Im\left\{\frac{\partial }{\partial t}\ln \label{transformfrequency} \left[W_\psi(t,s)\right]\right\}\\
\Upsilon_\psi(t,s)\equiv   \Re\left\{\frac{\partial }{\partial t}\ln \label{transformfrequency2} \left[W_\psi(t,s)\right]\right\}
\end{eqnarray}
which we refer to as the \emph{transform instantaneous frequency and bandwidth}.  We then have the estimates
\begin{eqnarray}
\widehat\omega^{\{\cdot\}}(t) \equiv \Omega_\psi\left(t,\widehat s^{\{\cdot\}}(t)\right)\\
\widehat\upsilon^{\{\cdot\}}(t) \equiv \Upsilon_\psi\left(t,\widehat s^{\{\cdot\}}(t)\right)
\end{eqnarray}
which are obtained from the values of the instantaneous frequency and bandwidth along the ridge curve.  In numeric implementation, differentiation is thus performed prior to the lookup along ridges, rather than the reverse.  The form of the instantaneous frequency estimate is found to be
\begin{multline}
\widehat\omega^{\{\cdot\}}(t)=\omega(t)\label{differentiatedestimate}
\left[1+P_\psi^2\left(\frac{1}{2} \frac{\omega''(t)}{\omega^ 3(t)}+\frac{\upsilon(t)}{\omega(t)}\frac{\omega'(t)}{\omega^ 2(t)}\right)\right.\\\left.+ O\left\{\epsilon_{\psi,3}^{\{\cdot,t\}}(t)\right\}\right]
\end{multline}
while that of the bandwidth is
\begin{multline}
\widehat\upsilon^{\{\cdot\}}(t)\label{bandwidthestimate}
=\upsilon(t)+\omega(t)\times \\\left[
P_\psi^2\left(\frac{1}{2}\frac{\upsilon''(t)}{\omega^3(t)}+
\frac{\upsilon(t)}{\omega(t)}\frac{\upsilon'(t)}{\omega^2(t)}\right) +O(\delta_{N_T}^4)+ O\left\{\epsilon_{\psi,3}^{\{\cdot,t\}}(t)\right\}\right]
\end{multline}
as follow at once from (\ref{imaglogtransform}) of Appendix~\ref{section:ridgecurves}.  The residual quantities $\epsilon_{\psi,3}^{\{\cdot,t\}}(t)$ are given by  (\ref{epsilontdef}).  The estimated instantaneous frequency and bandwidth plotted in Fig.~\ref{analytic-transforms}g--i have been constructed in this manner.

The leading-order perturbation term in the instantaneous frequency estimate (\ref{differentiatedestimate}) is at second order in $\delta_{N_T}$, but as $\upsilon(t)/\omega(t)$ is itself a first-order quantity, the estimated bandwidth in (\ref{bandwidthestimate}) is perturbed at first order in $\delta_{N_T}$.  The wavelet ridge estimation can therefore recover the instantaneous frequency with greater fidelity than it can the bandwidth.  The instantaneous frequency estimate (\ref{differentiatedestimate}) also has the desirable property of being identical for both the amplitude ridge curves and phase ridges curves, unlike the direct estimates of instantaneous frequency (\ref{directestimate}) which differ at order $\delta_{N_T}^2$.  Inserting (\ref{phasecurve}) into  (\ref{directestimate}) shows that the instantaneous frequency estimate (\ref{differentiatedestimate}) using either type of ridge curve is identical at leading order to that for direct estimate (\ref{directestimate}) of instantaneous frequency using the phase ridge curve.

For reference, it is useful to compare the estimated instantaneous frequency and bandwidth with the rates of change of the amplitude and phase of the localized analytic signal $x_\psi(t)$.  One finds
\begin{multline}
\omega_\psi(t)= \Im\left\{\frac{d}{dt}\ln \left[ x_\psi (t)\right]\right\}\\=\omega(t)\left[1 + \frac{1}{2}P_\psi^2  \frac{\omega''(t)}{\omega^3(t)} \right]+O(\delta_{N_T}^3)+O\left\{\epsilon_{\psi,3}^{\{t\}}\left(t,\omega_\psi/\omega(t)\right)\right\}
\label{differentiatelocally}
\end{multline}
for the rate of change of the phase and
\begin{multline}
\upsilon_\psi(t)\equiv \Re\left\{\frac{d}{dt}\ln \left[ x_\psi (t)\right]\right\}
=\upsilon(t)+\omega(t)\times \\\left[ P_\psi^2\left(\frac{1}{2}\frac{\upsilon''(t)}{\omega^3(t)}+
\frac{\upsilon(t)}{\omega(t)}\frac{\upsilon'(t)}{\omega^2(t)}\right)+O(\delta_{N_T}^4)\right.\\\left. + O\left\{\epsilon_{\psi,3}^{\{t\}}\left(t,\omega_\psi/\omega(t)\right)\right\}\right]
\end{multline}
for the relative rate of change of amplitude.  The form of the residual term is given in (\ref{partialtepsilon}).  Comparison with (\ref{differentiatedestimate}) and (\ref{bandwidthestimate}) shows that, while the instantaneous bandwidth estimate is identical to lowest perturbation order to the rate of change of amplitude of the localized analytic signal $x_\psi(t)$, the same is not true for the instantaneous frequency estimate.

The additional term in $\widehat\omega^{\{\cdot\}}(t)$ compared with $\omega_\psi(t)$ reflects the joint effect of contemporaneous amplitude and frequency modulation.  It does not occur in $\omega_\psi(t)$  since
\begin{multline}
 \omega_\psi(t)= \frac{d}{dt}\,\phi_\psi(t) =
\frac{d}{dt}\,\Im\ln W_\psi\left(t,\omega_\psi/\omega(t)\right)
=\\\Omega_\psi\left(t,\omega_\psi/\omega(t)\right)
-\omega_\psi\frac{\partial }{\partial s}\left.\Im\ln W_\psi(t,s)\right|_{t, \omega_\psi/\omega(t)}\frac{\omega'(t)}{\omega^2(t)}
\end{multline}
and evaluating the term proportional to $\omega'(t)$ from (\ref{ridgetheorem}), we find it cancels a similar term in $\Omega_\psi\left(t,\omega_\psi/\omega(t)\right)$, leading to
(\ref{differentiatelocally}).  The difference between $\widehat\omega^{\{\cdot\}}(t)$ and $\omega_\psi(t)$ is therefore attributed to a contribution to the rate of change of phase $\phi_\psi(t)$ due to the motion of the instantaneous frequency curve across scales at a fixed time.  The estimated instantaneous frequency thus seems anomalous---compared with the estimated amplitude, phase, and instantaneous bandwidth---in that it is not completely controlled at lowest perturbation order by the localized analytic signal.

Similarly we can estimate the second-order instantaneous modulation function $\widetilde \rho_2(t)$ by defining
\begin{equation}
   \widetilde P_2(t,s)\equiv\frac{\Upsilon^2_\psi(t,s)}{\Omega^2_\psi(t,s)} +\frac{\frac{\partial}{\partial t}\Upsilon_\psi(t,s)}{\Omega^2_\psi(t,s)} +i\frac{\frac{\partial}{\partial t}\Omega_\psi(t,s)}{\Omega^2_\psi(t,s)}
\end{equation}
which is evaluated along a ridge, leading to
\begin{equation}
\widetilde\rho_2^{\{\cdot\}}(t)
\equiv    \widetilde P_2\left(t,\widehat s^{\{\cdot\}}(t)\right).
\end{equation}
Proceeding as with the instantaneous frequency and bandwidth estimates, one can find
\begin{equation}
\widetilde\rho_2^{\{\cdot\}}(t)=
\widetilde\rho_2(t)\left[1 +O(\delta_{N_T})+ \cdots\right].
\end{equation}
where the leading-order term is again at order $O(\delta_{N_T})$, and with the ellipses denoting a twice-differentiated residual term following the development for the once-differentiated residual term in Appendix~\ref{appendix:bounding}. In the example in Fig.~\ref{analytic-transforms}, panels (j--l) show three versions of  $\widetilde\rho_2^{\{\cdot\}}(t)$ estimated in this way using three different wavelets.

\subsection{Application}\label{section:application}

The above development has shown that if the local stability level $\delta_{N_T}$ is known and is small compared to unity, one can choose a suitable wavelet such that the ridge-based estimates of the analytic signal $x_+(t)$, its amplitude $a_+(t)$ and phase $\phi_+(t)$, its instantaneous frequency $\omega(t)$ and bandwidth $\upsilon(t)$, and its second-order instantaneous modulation function $\widetilde\rho_2(t)$, are all close to their true values. A difficulty of course is that the local stability level $\delta_{N_T}$ is not generally known in practice.  One can estimate $\delta_{N_T}$, but as was seen in Fig.~\ref{analytic-transforms}, the estimated degree of smoothness of the signal depends upon the choice of wavelet.

Nevertheless, by presuming that a certain estimated signal is in fact correct, one can gain insight into what estimated signals are reasonable.  We apply this approach to the signal estimates shown in Fig.~\ref{analytic-transforms}.  There are two considerations: whether the estimated signal is sufficiently smooth, i.e. is characterized by a sufficiently small $\delta_{N_T}$, and secondly the wavelet is suitable for the stability level of the estimate it produces.  For simplicity we take the time period $T$ to be the entire record, and set the truncation level to $N=2$ so that we only need consider the second-order modulation function.


In Fig.~\ref{analytic-transforms}j--l, we have plotted
$\widetilde\rho_2(t)$ and drawn horizontal lines at $4/P_{\beta,\gamma}^4$.   Since the lowest-order stability condition, as discussed in Section~\ref{section:morsesuitabilty}, is $P_{\beta,\gamma}\le \sqrt{2/\delta_{N_T}}$ or equivalently $\delta_{N_T}^2\le 4/P_{\beta,\gamma}^4$, we should see $\widetilde\rho_2(t)$---which is of order $\delta_{N_T}^2$ by assumption---be bounded by these lines if the wavelet is suitable for the signal. One should keep in mind that the contamination of the signal by noise will tend to increase the roughness, so the estimated values of $\widetilde\rho_2(t)$ are anticipated to be somewhat too large.  From inspection, we see that extension of $\widetilde\rho_2(t)$ outside of these lines is minor for (j), occasional for (k), and---although it is difficult to see in this plot--- extensive for (l).  Since the great difference in the values of $\widetilde\rho_2(t)$ and $4/P_{\beta,\gamma}^4$ for the three estimated signals makes it difficult to compare them visually, we calculated some statistics to characterize their levels of variability.  The mean values of the ratio $|\widetilde\rho_2(t)|/(P_{\beta,\gamma}^4/4)$ are $0.54$,  $0.74$, and $2.10$ for (j--l), respectively, while the corresponding median values are $0.47$, $0.57$, and $1.74$.  Since the suitability criteria require that this quantity be smaller than unity, it appears that the wavelets used in the third column have a time duration that is too long, and this estimate would therefore be expected to be of a poor quality.

Thus the level of variability in the third estimated signal clearly exceeds that expected from the suitability conditions.  Let us say that the smoothest estimated signal, in Fig.~\ref{analytic-transforms}c, is in fact the true analytic signal to be estimated. The extensive excursions of $\widetilde\rho_2(t)$ outside the dotted lines in  Fig.~\ref{analytic-transforms}l means that the wavelet used in this column, Fig.~\ref{analytic-morsies}c, is not suitable to analyze this signal.  The ridge analysis using this wavelet has produced an estimated signal which it would not be able to recover accurately, an unacceptable result.   On the other hand, the horizontal lines correspond to values of $\delta_{N_T}$ of $0.44$, $0.22$, and $0.01$, respectively.  From Fig.~\ref{analytic-morsies}j we see that the variability of the estimated signal in the first column requires the relatively large choice of $\delta_{N_T}=0.44$.  Thus the estimated signal in the first column is quite rough, and indeed the rapid fluctuations in Fig.~\ref{analytic-morsies}g would lead one to suspect that this estimate is contaminated by noise.

Assuming that the estimated signal is the true signal, we can iterate the estimation procedure and ask which of the iterated estimates shows the least error.  We solve for the median and mean values of $|[\widehat x_+(t)-x_+(t)]/x_+(t)|^2$, in which each of the three estimates signals in Fig.~\ref{analytic-transforms}j--l plays the role of the true signal $x_+(t)$.  The mean values of the iterated deviations are $0.040$, $0.036$, and $0.057$, respectively, while the median values are $0.024$, $0.014$, and  $0.022$.  This means that the wavelet used in the middle column is able to recover the estimated signal it produces with the greatest degree of fidelity. Thus while the true signal remains unknown, we can say from a quantitative analysis that the estimate in Fig.~\ref{analytic-transforms}b is to be preferred.

It was shown in Section~\ref{section:morsesuitabilty} that for $1\le\gamma\le 6$, generalized Morse wavelet derivatives of all orders will satisfy the wavelet suitability criteria provided the lowest-order condition, $P_{\beta,\gamma}\le\sqrt{2/\delta{N_T}}$, is also satisfied.  This implies that a range of $\gamma$ values could be chosen for a fixed $P_{\beta,\gamma}$ and yield similar results.  To check this, we compute the wavelet ridge estimates for the wavelets shown in Fig.~\ref{analytic-morsies}f and Fig.~\ref{analytic-morsies}g, which like that in Fig.~\ref{analytic-morsies}b have
$P_{\beta,\gamma}=3$, but with $\gamma=1$ and $\gamma=6$ respectively.  The results (not shown) reveal both of the wavelet ridge estimates are very close to that for $\gamma=3$, as expected. This reflects the fact that the error terms due to higher-order wavelet derivatives have been successfully contained to higher perturbation order by the wavelet suitability criteria.

\subsection {Implications for Choice of Wavelet}\label{section:application2}

In this section we show how the ideas developed in this paper guide the choice of wavelet appropriate to the analysis of a given signal, using the generalized Morse wavelets as the reference point \cite{olhede02-itsp,lilly09-itsp}.

The higher-order wavelet properties will be addressed first.  The wavelet suitability conditions imply $1\le\gamma\le 6$ for the generalized Morse wavelets.  The $\gamma=3$ wavelets are in a sense optimal for fixed $P_{\beta,\gamma}$ since they have $\widetilde\Psi_{3;\beta,\gamma}(\omega_{\beta,\gamma})=0$ \cite{lilly09-itsp}. Thus after the leading-order terms proportional to $\widetilde\rho_2(t)$, the set of terms proportional to  $\widetilde\rho_3(t)$ vanishes, and the next contribution involves \emph{fourth-order} signal variability.  The second-order expansions we have emphasized are therefore particularly accurate for the $\gamma=3$ family.  By contrast, it is apparent that an inopportune choice of analyzing wavelet could lead to a very poor estimate of the analytic signal.  If one fixes $P_{\beta,\gamma}=\sqrt{\beta \gamma}$ and lets $\gamma$ increase without bound as $\beta$ decays to zero, violating the suitability conditions, (\ref{morsepsi3}) shows that the magnitude of $\widetilde\Psi_{3;\beta,\gamma}(\omega_{\beta,\gamma})$ increases without bound.  This implies there is no limit to how poor the estimated analytic signal can become if one chooses as an analyzing wavelet a function that, while analytic, is extremely asymmetric.

Nevertheless, the wavelet suitability criteria give some latitude in the choice of $\gamma$.  To understand why one might choose a particular value of $\gamma$ we consider first the roles of $\beta$ and $\gamma$ more generally.  We have noted that $\beta$ controls the time decay, with $\psi(t)/\psi(0)\sim|t|^{-(\beta+1)}$.  Meanwhile, $\gamma$ controls the high-frequency decay, as is clear from the frequency-domain form (\ref{morse}).  Thus increasing $\gamma$ is attractive if one would like to extend the analysis closer to the Nyquist frequency.  Decreasing $\gamma$ with $P_{\beta,\gamma}$ held fixed, on the other hand, lets $\beta$ and hence the rate of time decay be increased, minimizing leakage from distant times.  Further details on the roles of $\beta$ and $\gamma$ in setting the wavelet properties were considered by \cite{lilly09-itsp}.

With the higher-order errors made small by the constraint $1\le \gamma\le 6$, the leading-order error term is controlled by the choice of $P_{\beta,\gamma}$.  Our  analysis suggests that to minimize the leading-order error term, the wavelet should be chosen to be as short as possible---that is, $P_{\beta,\gamma}$ should be minimized.  In the example discussed in Section~\ref{section:application2} it was seen that the presence of noise compels us to choose $P_{\beta,\gamma}$ large enough to stabilize the transform against random fluctuations. Further examination of the impact of noise is  outside the scope of this paper.  But a second factor opposing the desire to make  $P_\psi$ small is that a wavelet cannot be made vanishingly short and still have attractive properties as a bandpass filter.

A conventional measure of the time-domain spread of a wavelet is its second moment
\begin{equation}
\sigma_{t;\beta,\gamma}^ 2 \equiv  \frac{\int_{-\infty}^\infty t^ 2|\psi_{\beta,\gamma}(t)|^ 2\,dt}{\int_{-\infty}^\infty|\psi_{\beta,\gamma}(t)|^ 2\,dt}.\label{timespread}
\end{equation}
Clearly the generalized Morse wavelets only have finite time spread $\sigma_{t;\beta,\gamma}^ 2$ for $\beta>1/2$, since $\beta=1/2$ implies $\psi(t)/\psi(0)\sim|t|^{-(3/2)}$, in which case the integrand in the numerator (\ref{timespread}) is proportional to $t^{-1}$.  Such long time decay is useless in practice.  At $\beta=1$, the time decay of the wavelet is already relatively slow at $t^{-2}$, and so this in some sense represents a lower bound for useful value of $\beta$.  Then the smallest value of $P_{\beta,\gamma}$ satisfying the wavelet suitability conditions would occur at $\gamma=1$, where we have $P_{\beta,\gamma}=1$.   This implies the wavelet executes one full cycle within its central window, an intuitive lower bound on the duration of a signal which is supposed to be a modulated oscillation.  As a result analysis of highly variable signals with $\delta_{N_T}$ of order unity will be problematic, but as mentioned earlier, such signals are not aptly described as modulated oscillations in the first place.

\section {Discussion}\label{conclusion}

This work has derived fundamental properties of the continuous analytic wavelet transform (AWT). In particular we have calculated an exact form for the AWT of a signal which may depart substantially from the case of negligible modulation. The key to achieving this representation is an expansion of the signal in terms of a set of appropriate time-varying functions---the \emph{instantaneous modulation functions}---which quantify the local degree of departure of the signal from a constant-amplitude, constant-frequency sinusoid.  The AWT is found to involve a series of interactions of increasingly higher-order instantaneous modulation functions of the signal with increasingly higher-order frequency-domain derivatives of the wavelet, a result termed the AWT representation theorem.  For signals or time intervals of a signal which are locally oscillatory, the AWT simplifies substantially.  By constraining the magnitude of frequency-domain derivatives of the wavelet, the Taylor expansion of the AWT with respect to scale can be reduced to a handful of important terms in the vicinity of the instantaneous frequency curve.

Wavelet ridge analysis, a means for estimating the properties of a modulated oscillation, was then revisited in the light of these results. Extending earlier work bounding the bias terms globally, we identified the lowest-order \emph{time-varying} bias of the estimated signal properties when the amplitude and frequency modulation are not negligible.  It was seen that amplitude- and phase-based ridge definitions are  different from one another, but that this difference is fact of secondary importance.  The leading-order error is instead due to the smoothing of the analytic signal by the wavelet along the signal's instantaneous frequency curve, an object we term the \emph{localized analytic signal}, not to the deviation of the instantaneous frequency curve from either type of ridge. In fact to leading perturbation order the estimated analytic signal is identical to the localized analytic signal. Amplitude and phase, as well as instantaneous bandwidth and frequency may be estimated with fidelity provided the signal modulation is not too strong and the wavelet is chosen appropriately.

Given the ubiquity of modulated oscillatory signals in a number of applications, these results will enable better characterization of such signals, and will add to the theory underpinning existing estimation methods.  For example, the discrete complex-valued decompositions mentioned in the Introduction have useful properties because they approximate a wavelet transform with an analytic mother wavelet function.  Our results are applicable to most of these decompositions, up to some (small) corrective error term which decreases with increasing scale or length of wavelet.  The Dual-Tree Complex Wavelet Transform (DCWT) \cite{selesnick05-ispm} is one such method. The DCWT has become an extremely popular tool in signal analysis, because it alleviates several shortcomings of the real discrete wavelet transform but maintains a controlled level of redundancy. Applications based on the observed properties of the DCWT coefficients are usually derived from its magnitude and phase properties \cite{selesnick05-ispm}. By applying our derived understanding of the AWT we may quantify aspects of the behavior of the DCWT. Another transform whose higher-order properties can be approximately determined from the results derived in this article is the chirplet transform \cite{mann95-itsp,candes02}. Thus although the primary goal of this paper is to understand and improve wavelet-based estimates of oscillatory signals, the results should also find applicability to other local estimation methods of modulated signals.

\appendices

\section{A Freely Distributed Software Package}\label{section:software}
All software associated with this paper is distributed as a part of a freely available Matlab toolbox called Jlab, written by the first author and available at \url{http://www.jmlilly.net}.  The Jsignal module of Jlab includes numerous routines for high-quality wavelet ridge analysis suitable for large data sets. Given an analytic signal, \texttt{instfreq} constructs the instantaneous frequency, bandwidth, and second-order modulation function.  The generalized Morse wavelets are implemented with \texttt{morsewave}, while their basic properties, peak frequency, and frequency-domain derivatives are computed in \texttt{morseprops}, \texttt{morsefreq}, and \texttt{morsederiv} respectively. The time spread of the generalized Morse wavelets is computed by \texttt{morsebox}, and \texttt{bellpoly} computes the Bell polynomials. The wavelet transform is implemented by \texttt{wavetrans} while \texttt{ridgewalk} has an efficient algorithm for finding the ridges.  All routines are well-commented and many have built-in automated tests or sample figures.  Finally, \texttt{makefigs\!\_\,analytic} generates all figures in this paper.

\section{Proof of the AWT Representation Theorem}\label{appendix:representation}

In this section we will use the notation $\psi_s(t)\equiv \psi(t/s)/s$ for a rescaled version of the wavelet.   Inserting the local modulation expansion (\ref{zBellN}) into the wavelet transform (\ref{wavetrans3}), one obtains
\begin{equation}
W_\psi(t,s)=W_{\Sigma_N}(t,s)+W_{R_{N+1}}(t,s)\label{transformpartition2}
\end{equation}
where [with $\widetilde\rho_0(t)\equiv 1$]
\begin{multline}
W_{\Sigma_N}(t,s)  \equiv \frac{1}{2} \,x_+(t)\times \\ \int_{-\infty}^{\infty} \psi_s^*\left( \tau\right) \,e^ {i\omega(t) \tau}
\sum_{n= 0}^N \frac{1}{n!}\left[\omega(t)\tau\right]^n\widetilde\rho_n(t)\,d\tau\label{WSigma}
\end{multline}
is the wavelet transform of the $N$th-order time-domain polynomial from the instantaneous modulation function signal expansion (\ref{zBellN}).  Note that in (\ref{transformpartition2}) the anti-analytic contribution vanishes on account of the analyticity of the wavelet.  $W_{R_{N+1}}(t,s)$  is implicitly defined by (\ref{transformpartition2}) as the difference between the wavelet transform of the signal $W_\psi(t,s)$ and the wavelet transform of the expansion  $W_{\Sigma_N}(t,s)$.

Now, the large-time decay of the wavelets is $O\left(t^{-r_\psi}\right)$ [see~(\ref{rdef})].  In order for the integrand in (\ref{WSigma}) to be square integrable, it is clear we must have $N\leq r_\psi-2$. Assuming that to be the case, (\ref{WSigma}) can simplify by substituting (\ref{nthderivative}) for the $n$th dimensionless derivative $\widetilde\Psi_n(\omega)$.  Then (\ref{WSigma}) becomes
\begin{equation}
W_{\Sigma_N}(t,s) =\frac{1}{2}\,x_+(t)\Psi^*(s\omega(t))
\sum_{n=0}^N \frac{(-i)^n  \widetilde\rho_n(t)}{n!} \,\widetilde{\Psi}_n^*(s\omega(t))\label{almostzBelltransform}
\end{equation}
and if we furthermore denote
\begin{equation}
\varepsilon_{\psi,N+1}(t,s)\equiv \frac{W_{R_{N+1}}(t,s)}{\frac{1}{2}x_+(t)\Psi^*\left(s\omega(t)\right)}\label{newvarepsilondef}
\end{equation}
the AWT representation theorem (\ref{zBelltransform}) follows by combining (\ref{transformpartition2}), (\ref{almostzBelltransform}), and (\ref{newvarepsilondef}).

To obtain bounds on the residual term $W_{R_{N+1}}(t,s)$, we split the wavelet transform integration into an inner and outer portion. Choose an energy level $\alpha$, with $1-\alpha<<1$, which determines a wavelet half-width $L_\psi(\alpha)$ as defined by (\ref{Ldef}).  We then write
\begin{multline}
W_{R_{N+1}}(t,s)=W_{I,R_{N+1}}(t,s;\alpha)+W_{O}(t,s;\alpha)\\-W_{O,\Sigma_{N}}(t,s;\alpha)
\end{multline}
where ``I'' and ``O'' denote integrations over the inner and outer ranges, respectively. The first of these three terms
\begin{multline}
W_{I,R_{N+1}}(t,s;\alpha)\\\equiv\frac{1}{2}\, x_+(t)  \int_{-sL_\psi(\alpha)}^{sL_\psi(\alpha)} \psi_s^*\left( \tau\right) \,e^ {i\omega(t) \tau } R_{N+1}(\tau,t)\,d\tau
\end{multline}
gives the integral of the residual term $ R_{N+1}(\tau,t)$, defined in (\ref{RNdef}), over the inner range. The second term
\begin{multline}
W_{O}(t,s;\alpha)\equiv \frac{1}{2} \int_{-\infty}^{-sL_\psi(\alpha)} \psi_s^*\left( \tau\right) x_+(t+\tau) \,d\tau\\+\frac{1}{2}
\int_{sL_\psi(\alpha)}^{\infty}\psi_s^*\left(\tau\right) x _+(t+\tau) \,d\tau
\end{multline}
is the integral including the entire signal over the outer range. The third term
\begin{multline}
W_{O,\Sigma_N}(t,s;\alpha)\equiv   \sum_{n= 0} I_{n}(t,s;\alpha) \equiv\\\frac{1}{2} \,x_+(t)\int_{-\infty}^{-sL_\psi(\alpha)} \psi_s^*\left( \tau\right) \,e^ {i\omega(t)\tau }\sum_{n= 0}^N \frac{1}{n!}\left[\omega(t)\tau\right]^n\widetilde\rho_n(t)\,d\tau\\+ \frac{1}{2} \,x_+(t)
\int_{sL_\psi(\alpha)}^{\infty}\psi_s^*\left( \tau\right) \,e^ {i\omega(t) \tau }\sum_{n= 0}^N \frac{1}{n!}\left[\omega(t)\tau\right]^n\widetilde\rho_n(t)\,d\tau
\label{outersummation}
\end{multline}
is the integral of the summation over the outer range; here we have also defined the contribution from the $n$th term in the summation $I_{n}(t,s;\alpha)$. The reason for this seemingly circuitous route to obtaining a bound is that we have only assumed derivatives of $x_+(t)$ to exist on the time interval $|t|\le sL_\psi(\alpha)$.  Outside this interval the residual $R_{N+1}(\tau,t)$ is no longer given by (\ref{RNdef}), but is still defined implicitly as the difference between the time series and the summation.

One finds the squared magnitude of the inner term is subject to the bound
\begin{multline}
\left|W_{I,R_{N+1}}(t,s;\alpha)\right|^2\\\le \frac{1}{4} \, \left|x_+(t)\right|^2 \alpha\,\frac{c_\psi^2}{s}
\sup_{\tau \in \left[-s L_\psi(\alpha),s L_\psi(\alpha)\right]}
\left|R_{N+1}(\tau,t)\right|^ 2
\end{multline}
which increases with increasing $\alpha$; here $c_\psi^2$  is the wavelet energy
\begin{equation}
c_\psi^2\equiv \int_{-\infty}^{\infty} \left|\psi(t)\right|^2\,dt=\frac{1}{2\pi}\int_{-\infty}^{\infty} \left|\Psi(\omega)\right|^2\,d\omega.
\end{equation}
By the triangle inequality together with the Cauchy-Schwarz inequality we also find that
\begin{eqnarray}
\left|W_{O}(t,s;\alpha)\right|^2&\le& \frac{1}{4}\,\|x_+\|^2\left(\frac{1-\alpha}{2}\right)\frac{c_\psi^2}{s}
\end{eqnarray}
[with $\|x_+\|$ denoting the $L^2$ norm of $x_+(t)$] and thus the contribution to the wavelet transform from $|t|> sL_{\psi}(\alpha)$ is negligible if $\alpha$ is chosen to be sufficiently close to unity. The contributions of these two terms are therefore antagonistic.

To find the bound on the third term, we define $b_\psi>0$ and $d_\psi>0$ to be constants chosen such that
\begin{eqnarray}
\left|\psi(t)\right|&\le& b_\psi |t|^{-r_\psi}\\
\left|\psi(t)\right|&\sim& d_\psi |t|^{-r_\psi}
\end{eqnarray}
where $r_\psi$ gives the wavelet time decay.  Then one may note
\begin{equation}
L_\psi^{-1}(\alpha)\approx\left(\frac{1}{2}(1-\alpha)(2r_\psi-1)\frac{c_\psi^2}{d_\psi^2}\right)^
{1/(2r_\psi -1)}
\end{equation}
and from this we find
\begin{multline}
\left|I_{n}(t,s;\alpha)\right|^ 2\leq  \frac{1}{2}\left|x_+(t)\right|^ 2\frac{b_\psi^ 2}{s}\left|\frac{[s\omega(t)]^n}{n!}\right|^ 2\left|\widetilde\rho_n(t)\right|^ 2\times \\\frac{\left(L_\psi^{-1}(\alpha)\right)^{2(r_\psi-n)-1}}{2(r_\psi-n) -1}
\end{multline}
which follows in a few lines of algebra from (\ref{outersummation}) using the triangle inequality together with the observation and also (\ref{alphadef}) and (\ref{rdef}).  Then finally
\begin{equation}
\left|W_{O,\Sigma_N}(t,s;\alpha)\right|^ 2\leq \sum_{n= 0}^N \left|I_{n}(t,s;\alpha)\right|^ 2
\end{equation}
by the triangle inequality. The three components of $W_{R_{N+1}}(t,s)$ are therefore bounded, and $W_{R_{N+1}}(t,s)$ itself is bounded by another application of the triangle inequality.

\section{Proof of the AWT Scale Deviation Expansion}\label{section:deviation}
Noting then the normalized wavelet derivatives have the Taylor series expansion
\begin{eqnarray}
\widetilde \Psi_{n}(s\omega)& = &  \sum_{m= 0}^\infty \frac{1}{m!} \widetilde \Psi_{n+m}(\omega_\psi) \left(\frac{s\omega}{\omega_\psi}-1\right)^m\label{normalizedexpansion}
\end{eqnarray}
we insert this into the wavelet transform representation theorem given in (\ref{zBelltransform}), yielding 
\begin{multline}
 W_\psi(t,s) =
x_+(t)
\left[\sum_{n= 0}^N\sum_{m= 0}^\infty  \frac{(-i)^n}{n! m!}\,\widetilde\Psi^*_{m+n}(\omega_\psi)\,\widetilde\rho_n(t) \times  \right.\\\left.
\left(\frac{s\omega(t)}{\omega_\psi}\right)^n \label{doubleexpansion} \left(\frac{s\omega(t)}{\omega_\psi}-1\right)^m+ \varepsilon_{\psi,N+1}(t,s)\right]
\end{multline}
where we let $\widetilde \Psi_ 0(\omega_\psi)\equiv  1$ and $\widetilde\rho_ 0(t)\equiv  1$. Now expanding powers of $s\omega(t)/\omega_\psi$  via the binomial theorem, one finds
\begin{multline}
\left(\frac{s\omega(t)}{\omega_\psi}\right)^n =\left( \frac{s\omega(t)}{\omega_\psi}-1+1\right)^n 
\\=\sum_{p= 0}^n \frac{n!}{(n-p)!p!}\left(\frac{s\omega(t)}{\omega_\psi}-1\right)^p
\end{multline}
and inserting into (\ref{doubleexpansion}),  we obtain the triple summation (\ref{triplesummation}).  Writing out all terms up to $m=2$, $n=2$ leads to
\begin{multline}
\frac{ W_\psi(t,s)}{x_+(t)}=1+\frac{1}{2}\widetilde\Psi^*_{2}(\omega_\psi)\left[\Delta\omega(t,s)\right]^2
+\ldots \\
-i\widetilde\rho_1(t) \label{allwritten} \left[\widetilde\Psi^*_{2}(\omega_\psi)\Delta\omega(t,s)\hspace{1.2in}\right.\\\left.
+\left(\widetilde\Psi^*_{2}(\omega_\psi)+\frac{1}{2}\widetilde\Psi^*_{3}(\omega_\psi)\right)\left[\Delta\omega(t,s)\right]^2+\ldots\right]\\
- \frac{1}{2} \widetilde\rho_2(t) \left[\widetilde\Psi^*_{2}(\omega_\psi)+
\left(2\widetilde\Psi^*_{2}(\omega_\psi)+\widetilde\Psi^*_{3}(\omega_\psi)\right)\Delta\omega(t,s)
\right.\\\left.+\left(\widetilde\Psi^*_{2}(\omega_\psi)+2\widetilde\Psi^*_{3}(\omega_\psi)+\frac{1}{2}\widetilde\Psi^*_{4}(\omega_\psi)\right)\left[\Delta\omega(t,s)\right]^2
+\ldots\right]\\+\varepsilon_{\psi,3}(t,s)
\end{multline}
where $\varepsilon_{\psi,3}(t,s)$  captures the influence of terms of higher order in the instantaneous modulation functions  $\widetilde\rho_n(t)$, while the ellipses denote terms of higher order in the scale deviation $\Delta\omega(t,s)\equiv s\omega(t)/\omega_\psi-1$.  Note that we have used the facts that the first derivative of the wavelet at the peak frequency $\widetilde\Psi_{1}(\omega_\psi)$ vanishes by definition, and that  $\Psi(\omega_\psi)\equiv 2$.

\section{Bounding the Differentiated Residual}\label{appendix:bounding}

For later use it will be necessary to obtain expansions for time and scale derivatives of the transform residual term $\varepsilon_{\psi,N+1}(t,s)$ (\ref{newvarepsilondef}), which we denote by
\begin{eqnarray}
\varepsilon_{\psi,N+1}^{\{t\}}(t,s) &\equiv &\label{partialtepsilon0}
 \frac{s}{\omega_\psi} \frac{\partial}{\partial t}\, \varepsilon_{\psi,N+1}(t,s)\\
\varepsilon_{\psi,N+1}^{\{s\}}(t,s) &\equiv &\label{partialsepsilon0}
 s\frac{\partial}{\partial s}\, \varepsilon_{\psi,N+1}(t,s).
\end{eqnarray}
Note that the additional factor of $\omega_\psi$ in the denominator of the former is convenient since it renders the derivative dimensionless, like the scale derivative (since $s$ is itself dimensionless).

Throughout this section we assume $\left(t,s\right)\in \calR_{2;\psi}(\delta_{N_T} )$ so that $s\omega(t)/\omega_\psi=1+O\left(\delta_{N_T}^2\right)$.  The time derivative of the transform residual is
\begin{multline}\label{longvarexpression}
 \varepsilon_{\psi,N+1}^{\{t\}}(t,s)=
-  \varepsilon_{\psi,N+1}(t,s)\times \\ \frac{s}{\omega_\psi} \,\left\{\upsilon(t)+i\omega(t)+\widetilde \Psi_1(s\omega(t)) \frac{\partial}{\partial t}\ln\left[\omega(t)\right]\right\}\\+
\frac{ \frac{s}{\omega_\psi}  \frac{\partial}{\partial t} \left\{W_{R_{N+1}}(t,s)\right\}}{\frac{1}{2}\,x_+(t)\Psi^*(s\omega(t))}
\end{multline}
but from (\ref{normalizedexpansion}) one has
\begin{equation}
\widetilde \Psi_1(s\omega(t)) = O\left\{\widetilde \Psi_2(\omega_\psi)\times \left(\frac{s\omega(t)}{\omega_\psi}-1\right)\right\}
 = O\left(\delta_{N_T}\right)\label{psi1thingy}
\end{equation}
for $\left(t,s\right)\in \calR_{2;\psi}(\delta_{N_T} )$. Recalling also (\ref{smalldrho}), the term on the second line in (\ref{longvarexpression}) is seen to be an order unity quantity, and we then find
\begin{multline}
\varepsilon_{\psi,N+1}^{\{t\}}(t,s)=O\left(\varepsilon_{\psi,N+1}(t,s)\right)+
\frac{\frac{s}{\omega_\psi} \frac{\partial}{\partial t} \, \left\{W_{R_{N+1}}(t,s)\right\}}{\frac{1}{2}\,x_+(t)\Psi^*(s\omega(t))}
\end{multline}
where the order of the second term remains to be found.  Likewise for the scale derivative one obtains
\begin{multline}
\varepsilon_{\psi,N+1}^{\{s\}}(t,s)=\\- s \widetilde\Psi_1(s\omega(t))  \times \varepsilon_{\psi,N+1}(t,s)\label{scalederivative}
+\frac{s\frac{\partial}{\partial s}\left\{W_{R_{N+1}}(t,s)\right\}}{\frac{1}{2}\,x_+(t)\Psi^*(s\omega(t))}
\end{multline}
but on account of (\ref{psi1thingy}), this becomes
\begin{multline}
\varepsilon_{\psi,N+1}^{\{s\}}(t,s)=\\O\left(\delta_{N_T} \times \varepsilon_{\psi,N+1}(t,s)\right)
+\frac{s\frac{\partial}{\partial s}\left\{W_{R_{N+1}}(t,s)\right\}}{\frac{1}{2}\,x_+(t)\Psi^*(s\omega(t))}
\end{multline}
again leaving the order of the second term to be found.

We can obtain bounds for the numerators in the preceding expressions as follows. Define
\begin{eqnarray}
U_\psi(t,s)& \equiv & \frac{s}{\omega_\psi} \frac{\partial}{\partial t}\, W_\psi(t,s)\\
V_\psi(t,s)& \equiv & s \frac{\partial}{\partial s}\, W_\psi(t,s)
\end{eqnarray}
and note that $U_\psi(t,s)$ and $V_\psi(t,s)$ may themselves be written as wavelet transforms using modified wavelets. Define
\begin{eqnarray}
\theta(t)& \equiv &-\psi'(t)/\omega_\psi\\
\varphi(t)& \equiv &-\left[\psi(t)+t\psi'(t)\right]\label{varwavelet}
\end{eqnarray}
having Fourier transforms $\Theta(\omega)=-i(\omega/\omega_\psi)\Psi(\omega)$ and $\Phi(\omega)=\omega\frac{d}{d\omega}\Psi(\omega)$ respectively. The differentiated wavelet transforms may then be written
\begin{eqnarray}
U_\psi(t,s)& =& \int_{-\infty}^{\infty} \frac{1}{s} \theta^*\left( \frac{\tau-t}{s}\right) x(\tau)\,d\tau \label{wavetransU}\\
V_\psi(t,s)& = & \int_{-\infty}^{\infty} \frac{1}{s} \varphi^*\left( \frac{\tau-t}{s}\right) x(\tau)\,d\tau \label{wavetransV}
\end{eqnarray}
using the definition of the wavelet transform (\ref{wavetrans}).  Note that by incorporating the derivatives into the wavelets, the original signal remains in both integrands.

The functions $\theta(t)$ and $\varphi(t)$ are valid wavelets provided that they have finite energy and that the Fourier transform $\Psi(\omega)$ of the original wavelet satisfies
\begin{eqnarray}
\int_{-\infty}^{\infty}\left|\omega\right| \left|\Psi(\omega)\right|^ 2\,d\omega &<&\infty\label{uadmiss}\\
\int_{-\infty}^{\infty}\left|\omega\right| \left|\frac{d}{d\omega}\Psi(\omega)\right|^ 2\,d\omega &<&\infty\label{vadmiss}
\end{eqnarray}
which together constitute the admissibility conditions for $\theta(t)$ and $\varphi(t)$ respectively. Now inserting the local modulation expansion of the signal (\ref{zBellN}) into (\ref{wavetransU}) and (\ref{wavetransV}) we obtain [mirroring the development of Appendix~\ref{appendix:representation}]
\begin{eqnarray}
\frac{s}{\omega_\psi}\frac{\partial}{\partial t}\, W_\psi(t,s)&=&U_{\Sigma_N}(t,s)+U_{R_{N+1}}(t,s)\label{transformpartitionU}\\
s \frac{\partial}{\partial s}\, W_\psi(t,s)&=&V_{\Sigma_N}(t,s)+V_{R_{N+1}}(t,s)\label{transformpartitionV}
\end{eqnarray}
where the individual terms are defined analogously to those in the original transform of the signal as in (\ref{transformpartition2}).  In order for the integrals implied on the right-hand side to be well defined, we must have the truncation level $N$ satisfy $N\leq r_\theta-2$ and $N\leq r_\varphi-2$, where $r_\theta$ and $r_\varphi$ are the long-time decay of the differentiated wavelets defined as in (\ref{rdef}).  Henceforth we assume this to be the case.

Since by construction, $U_{\Sigma_N}(t,s)$ on the right-hand side of (\ref{transformpartitionU}) is equal to the derivative of the summation term that is implicit on the left-hand side, and similarly for $V_{\Sigma_N}(t,s)$ in (\ref{transformpartitionV}), we may also note
\begin{eqnarray}
 \frac{s}{\omega_\psi}\frac{\partial}{\partial t}\, W_{R_{N+1}}(t,s)&=&U_{R_{N+1}}(t,s)\label{transformpartitionU2}\\
s \frac{\partial}{\partial s}\, W_{R_{N+1}}(t,s)&=&V_{R_{N+1}}(t,s)\label{transformpartitionV2}.
\end{eqnarray}
The differentiated residuals (\ref{partialtepsilon0}) and (\ref{partialsepsilon0})  thus become
\begin{equation}
\varepsilon_{\psi,N+1}^{\{t\}}(t,s)=\label{partialtepsilon}
 O\left(\varepsilon_{\psi,N+1}(t,s)\right)+
\frac{ U_{R_{N+1}}(t,s)}{\frac{1}{2}\,x_+(t)\Psi^*(s\omega(t))}
\end{equation}
for the time differentiation and
\begin{equation}\varepsilon_{\psi,N+1}^{\{s\}}(t,s)=
 O\left(\varepsilon_{\psi,N+1}(t,s)\right)+\label{partialsepsilon}
\frac{V_{R_{N+1}}(t,s)}{\frac{1}{2}\,x_+(t)\Psi^*(s\omega(t))}
\end{equation}
for the scale differentiation.  $U_{R_{N+1}}(t,s)$ and $V_{R_{N+1}}(t,s)$ may then be bounded in the same manner as for $W_{R_{N+1}}(t,s)$ in Appendix~\ref{appendix:representation}, but using the modified wavelets $\theta(t)$ and $\varphi(t)$.

\section{Proofs of the Forms of the Ridge Curves}~\label{section:ridgecurves}
To obtain expressions for the ridge curves, it is necessary to assume at the outset the order of the deviation of a ridge from an instantaneous frequency curve.  We assume $\left(t,s^{\{\cdot\}}\right)\in \calR_{2;\psi}(\delta_{N_T} )$, i.e. that the ridge curve lies in the $2$-neighborhood of the instantaneous frequency curve; it is found that the ridge equations do indeed have solutions within this neighborhood.  Also, for convenience we take  $\widetilde\Psi_{n}(\omega_\psi)$ to be  real-valued for $n\leq 4$. For the amplitude ridges, we wish to solve (\ref{ampridge}--\ref{ampridge2}), while the phase ridges satisfy (\ref{phaseridge}--\ref{phaseridge2}).

Amplitude ridges will be considered first. Since $\ln(1+x)=x -x^2/2 +\dots$, we may write the log of the analytic wavelet transform as [from (\ref{transformpartition})]
\begin{multline}\label{partionexpansion}
\ln W_\psi(t, s)=\ln x_+(t) + \Delta x_\psi(t) + \Delta W_\psi(t, s) \\-\frac{1}{2} \left[ \Delta x_\psi(t) + \Delta W_\psi(t, s)\right]^2+\ldots
\end{multline}
where we note that the squared term cannot be neglected in what follows. Inserting  (\ref{locallyanalytic})  and  (\ref{ridgetheorem}) for $\Delta x_\psi(t)$ and $\Delta W_\psi(t,s)$, respectively, one obtains for the real part
\begin{multline}
\Re\left\{\ln W_\psi(t,s)\right\}=
 \Re\left\{\ln x_+(t)+\Delta x_\psi (t)-\frac{1}{2}\left[\Delta x_\psi (t)\right]^2\right\}\\-\Delta\omega(t,s)
 \left[ \Re\left\{\widetilde\rho_2(t)\right\}\left(\widetilde\Psi_{2}(\omega_\psi) +\frac{1}{2}\widetilde\Psi_{3}(\omega_\psi)\right) \right.\\\left. +\frac{1}{6}\Im\left\{\widetilde\rho_3(t)\right\}\widetilde\Psi_{4}(\omega_\psi)
  -\frac{1}{2}\Im\left\{\widetilde\rho_1(t)\widetilde\rho_2(t)\right\}\widetilde\Psi_{2}^ 2(\omega_\psi)
 \right] \\+\frac{1}{2}\left[\Delta\omega(t,s)\right]^ 2\widetilde\Psi_{2}(\omega_\psi)+ O\left(\delta_{N_T} ^{4}\right)+O\left(\varepsilon_{\psi,4}(t,s)\right)
\label{reallogtransform}
\end{multline}
with truncation level $N_T=3$.

In the 2-neighborhood  $(t,s)\in {\cal R}_{2;\psi}(\delta_{N_T})$ of an instantaneous frequency curve
\begin{equation}
s\frac{\partial}{\partial s}\Delta\omega(s,t)=\frac{s\omega(t)}{\omega_\psi}=O (1)
\end{equation}
so that taking a scale derivative of such terms transforms an $O(\delta_{N_T}^2)$ term into an $O(1)$ term.  Applying the scale derivative $s\frac{\partial}{\partial s}$ to (\ref{reallogtransform}) and evaluating the result along the ridge then leads to
\begin{multline}
- \left[\frac{\upsilon^2(t)}{\omega^2(t)} \left(\widetilde\Psi_{2}(\omega_\psi)+\frac{1}{2}\widetilde\Psi_{3}(\omega_\psi)\right) +\frac{1}{6}\Im\left\{\widetilde\rho_3(t)\right\}\widetilde\Psi_{4}(\omega_\psi)\right.\\\left.-\frac{1}{2}\frac{\upsilon(t)}{\omega(t)} \frac{\omega'(t)}{\omega^2(t)}  \widetilde\Psi_{2}^2(\omega_\psi)\right] +\left(\frac{\widehat s^{\,\{a\}}(t)\omega(t)}{\omega_\psi}-1\right)\widetilde\Psi_{2}(\omega_\psi)
\\+ O\left(\delta_{N_T} ^{2}\right)+O\left\{\delta_{N_T}\times\epsilon_{\psi,3}^{\{a,s\}}(t)\right\}= 0
\end{multline}
and dividing through by  $\widetilde\Psi_{2}(\omega_\psi)=O\left(\delta_{N_T} ^{-1 }\right)$  one obtains (\ref{amplitudecurve}) for the amplitude ridges.  The residual quantity in the above is defined as
\begin{equation}
\delta_{N_T}\times\epsilon_{\psi,N}^{\{\cdot,s\}}(t) \equiv  \varepsilon_{\psi,N+1}^{\{s\}}\left(t,s^{\{\cdot\}}\right)  \label{epsilonsdef}
\end{equation}
where an expression for $\varepsilon_{\psi,N+1}^{\{s\}}\left(t,s\right)$ is given by (\ref{partialsepsilon}) of Appendix~\ref{appendix:bounding}; here again the superscript ``$\cdot$'' could be either ``$a$'' or ``$p$''.

For the phase ridges, we proceed by defining the complex-valued transform quantity
\begin{equation}
H_{\psi}(t,s) \equiv \Omega_\psi(t,s)-i\Upsilon_\psi(t,s)\equiv
-i\frac{\partial}{\partial t} \ln W_\psi(t, s)\label{transformcomplexfrequency}
\end{equation}
in analogy with the signal's complex instantaneous frequency $\eta(t)\equiv\omega(t)-i\upsilon(t)$.  We then differentiate the wavelet transform $W_\psi(t,s)$ as expressed in (\ref{transformpartition}), including terms from $\Delta x_\psi(t)$ and as well as from $\Delta W_\psi(t,s)$.  The orders of the various terms can be assessed by recalling (\ref{smalldrho}) for derivatives of the instantaneous modulation functions, and by noting that for $(t,s)\in {\cal R}_{2;\psi}(\delta_{N_T})$ we have
\begin{multline}
\frac{\partial}{\partial t}\Delta\omega(t,s)
=\frac{s\omega'(t)}{\omega_\psi} = \omega(t)\frac{\omega'(t)}{\omega^2(t)}\left[1+O(\delta_{N_T}^2)\right]
\\=\omega(t)\times O(\delta_{N_T}^2)
\end{multline}
for time derivatives of the scale deviation.  With a truncation level of $N_T=2$, (\ref{transformcomplexfrequency}) becomes for $(t,s)\in {\cal R}_{2;\psi}(\delta_{N_T})$, by differentiating (\ref{partionexpansion}),
\begin{multline}
H_{\psi}(t,s) =\eta(t)\label{transformeta}
-\widetilde\Psi_{2}(\omega_\psi)\frac{\partial}{\partial t}\left[\Delta\omega(t,s)\widetilde\rho_1(t)-i\frac{1}{2}\widetilde\rho_2(t)\right]
\\+ \omega(t)\times O\left(\delta_{N_T}^3\right) +\omega(t)\times O\left(\varepsilon_{\psi,3}^{\{t\}}(t,s) \right)
\end{multline}
where the residual term is defined by (\ref{partialtepsilon}) of Appendix~\ref{appendix:bounding}.  Writing out terms we find
\begin{multline}
H_{\psi}\left(t,s\right) =\eta(t)-\omega(t)\times \\
\widetilde\Psi_{2}(\omega_\psi)\left[\frac{1}{2}\frac{\omega''(t)}{\omega^3(t)}+\frac{\upsilon(t)}{\omega(t)}\frac{\omega'(t)}{\omega^2(t)}-i\frac{1}{2}\frac{\upsilon''(t)}{\omega^3(t)}-i\frac{\upsilon(t)\upsilon'(t)}{\omega^3(t)}\right]
\\
+\omega(t)\times O\left(\delta_{N_T}^3\right)+\omega(t)\times O\left(\varepsilon_{\psi,3}^{\{t\}}(t,s) \right) \label{imaglogtransform}
\end{multline}
where the entire term in brackets is of second order in $\delta_{N_T}$.

Rearranging the phase ridge condition (\ref{phaseridge}) yields
\begin{equation}
\frac{s^{\,\{p\}}(t)\omega(t)}{\omega_\psi}\frac{1}{\omega(t)} \Omega_\psi\left(t,\widehat s^{\,\{p\}}(t)\right) =1
\end{equation}
and inserting the imaginary part of (\ref{imaglogtransform}), one finds
\begin{multline}
\frac{s^{\,\{p\}}(t)\omega(t)}{\omega_\psi}\left[1 -\left(\frac{1}{2}\frac{\omega''(t)}{\omega^3(t)} +\frac{\omega'(t)}{\omega^2(t)} \frac{\upsilon(t)}{\omega(t)} \right)\widetilde\Psi_{2}(\omega_\psi)\right.\\\left.+ O\left(\delta_{N_T} ^{3}\right)+O\left\{\epsilon_{\psi,3}^{\{p,t\}}(t)\right\}
\right]=1
\end{multline}
from which the form of the phase ridge curve (\ref{phasecurve}) follows.  Here we have defined
\begin{equation}
\epsilon_{\psi,N+1}^{\{\cdot,t\}}(t) \equiv  \label{epsilontdef}\varepsilon_{\psi,N+1}^{\{t\}}\left(t,s^{\{\cdot\}}\right)
\end{equation}
as the residual term along the ridge.

\section*{Acknowledgment}
We thank three anonymous reviewers for their constructive feedback.


\begin{biography}{Jonathan M. Lilly}
(M'05) was born in Lansing, MI,
in 1972. He received the B.S. degree in geology and
geophysics from Yale University, New Haven, CT,
in 1994, and the M.S. and Ph.D. degrees in physical
oceanography from the University of Washington,
Seattle, in 1997 and 2002, respectively.\\\indent
He was a Postdoctoral Researcher with the
Applied Physics Laboratory and School of Oceanography,
University of Washington, from 2002 to
2003, and with the Laboratoire d'Oc\'eanographie
Dynamique et de Climatologie, Universit\'e Pierre
et Marie Curie, Paris, France, from 2003 to 2005. Since 2005, he has been
a Research Associate with Earth and Space Research, a nonprofit scientific
institute in Seattle. His research interests are oceanic vortex structures, satellite
oceanography, time/frequency analysis methods, and wave--wave interactions.
\\\indent Dr. Lilly is a member of the American Meteorological Society and of the
American Geophysical Union.
\end{biography}

\begin{biography}{Sofia C. Olhede}
(M'06) was born in Spanga, Sweden, in 1977. She received the M. Sci. and
Ph.D. degrees in mathematics from Imperial College
London, London, U.K., in 2000 and 2003,
respectively.\\\indent
She was a Lecturer (2002--2006) and
Senior Lecturer (2006--2007) with the Mathematics
Department, Imperial College London. In 2007,
she joined the Department of Statistical Science, University
College London, where she is Professor of
statistics.
Her research interests include the analysis of complex-valued stochastic processes, nonstationary time series, and inhomogeneous random fields.
She is an Associate Editor of the \emph{Journal of the Royal Statistical
Society, Series B (Statistical Methodology)}.\\\indent
Prof. Olhede serves on the Research Section of the Royal Statistical
Society.
\end{biography}

\label{lastpage}
\end{document}